\documentstyle[12pt]{article}
\input{amssym.def}
\input{amssym.tex}
\parskip=\medskipamount
\newtheorem{definition}{Definition}[section]
\newtheorem{lemma}[definition]{Lemma}
\newtheorem{theorem}[definition]{Theorem}
\newtheorem{proposition}[definition]{Proposition}
\newtheorem{corollary}[definition]{Corollary}
\newtheorem{remark}[definition]{Remark}
\newtheorem{example}[definition]{Example}
\newtheorem{examples}[definition]{Examples}

 1   
\font\ddpp=msbm10  scaled \magstep 1  

\def\QED{\hskip0.1em\hfill\null\ \null\nobreak\hfill
\kern3pt\lower1.8pt\vbox{\hrule\hbox
{\vrule\kern1pt\vbox{\kern1.7pt \hbox{$\scriptstyle
QED$}\kern0.2pt}\kern1pt\vrule}\hrule}}
\def\R{\hbox{\ddpp R}}               

\def\H{\hbox{\ddpp H}}
\def\lcf{\lbrack\! \lbrack}
\def\rcf{\rbrack\! \rbrack}
\def\theequation{\arabic{section}.\arabic{equation}}

\newcommand\prueba {\mbox{{\em Proof: }}}

\begin{document}
\baselineskip=.55cm
\title{{\bf GENERALIZED LIE BIALGEBRAS AND JACOBI STRUCTURES ON LIE GROUPS}}
\author{David Iglesias, Juan C. Marrero
\\ {\small\it Departamento de Matem\'atica
Fundamental, Facultad de Matem\'aticas,}\\[-8pt] {\small\it
Universidad de la Laguna, La Laguna,} \\[-8pt] {\small\it
Tenerife, Canary Islands, SPAIN,}\\[-8pt] {\small\it E-mail:
diglesia@ull.es, jcmarrer@ull.es} }
\date{}

\maketitle
\baselineskip=.3cm
\begin{abstract}
{\small
We study generalized Lie bialgebroids over a single point, that is,
generalized Lie bialgebras. Lie bialgebras are examples of
generalized Lie bialgebras. Moreover, we prove that the last ones can
be considered as the infinitesimal invariants of Lie groups endowed
with a certain type of Jacobi structures. We also propose a method to
obtain generalized Lie bialgebras. It is a generalization of the
Yang-Baxter equation method. Finally, we describe the structure of a
compact generalized Lie bialgebra.
}
\end{abstract}
\begin{quote}
{\it Mathematics Subject Classification} (2000): 17B62, 22Exx, 53D05,
53D10, 53D17. 

{\it Key words and phrases}: Jacobi manifolds, Poisson manifolds,
contact structures, locally conformal symplectic structures, Lie
algebras, Lie bialgebras, Yang-Baxter equation, Lie groups.  
\end{quote}

\section{Introduction}
\baselineskip=.55cm
A Jacobi structure on a manifold $M$ is a 2-vector $\Lambda$ and a
vector field $E$ on $M$ such that $[\Lambda ,\Lambda ]=2E\wedge
\Lambda$ and $[E,\Lambda ]=0$, where $[\, ,\, ]$ is the
Schouten-Nijenhuis bracket \cite{Li2}. If $(M,\Lambda ,E)$ is a
Jacobi manifold one can define a bracket of functions, the Jacobi
bracket, in such a way that the space $C^\infty (M,\R )$ endowed with
the Jacobi bracket is a local Lie algebra in the sense of Kirillov
\cite{K}. Conversely, a local Lie algebra structure on $C^\infty
(M,\R )$ induces a Jacobi structure on $M$ \cite{GL,K}. Jacobi
manifolds are natural generalizations of Poisson, contact and locally
conformal symplectic manifolds.

A category with close relations to Jacobi geometry is that of Lie
algebroids. In fact, if $M$ is an arbitrary manifold, the vector
bundle $TM\times \R\to M$ possesses a natural Lie algebroid structure
and, moreover, if $M$ is a Jacobi manifold then the 1-jet bundle
$T^\ast M\times \R\to M$ admits a Lie algebroid structure \cite{KS}.
However, as Vaisman proved in \cite{V2}, the pair $(TM\times
\R,T^\ast M\times \R )$ is not a Lie bialgebroid in the sense of
Mackenzie and Xu \cite{MX} (or Kosmann-Schwarzbach \cite{K-S}). This
is an important difference with the Poisson case. Indeed, if $M$ is a
Poisson manifold, the vector bundle $T^\ast M\to M$ is a Lie
algebroid \cite{BV,CDW,F,V} and, in addition, if on the dual bundle
$TM\to M$ we consider the natural Lie algebroid structure then the
pair $(TM,T^\ast M)$ is a Lie bialgebroid \cite{MX}.

The above results (about the relation between Jacobi structures and
Lie bialgebroids) and some suitable examples of linear Jacobi
structures on vector bundles obtained in \cite{IM0} motivated the
introduction, in \cite{IM}, of the definition of a generalized Lie
bialgebroid, a generalization of the notion of a Lie bialgebroid. A
generalized Lie bialgebroid is a pair $((A,\phi _0),(A^\ast ,X_0))$,
where $A$ is a Lie algebroid over $M$, $\phi _0$ is a 1-cocycle in
the Lie algebroid cohomology complex of $A$ with trivial
coefficients, $A^\ast$ is the dual bundle to $A$ which admits a Lie
algebroid structure and $X_0$ is a 1-cocycle of $A^\ast$. Moreover,
the Lie algebroids $A$ and $A^\ast$ and the 1-cocycles $\phi _0$ and
$X_0$ must satisfy some compatibility conditions. When $\phi _0$ and
$X_0$ are zero, we recover the notion of a Lie bialgebroid. In
addition, if $(M,\Lambda ,E)$ is a Jacobi manifold we proved that the
pair $((TM\times \R,\phi _0),(T^\ast M\times \R,X_0))$ is a
generalized Lie bialgebroid, where $\phi _0=(0,1)\in \Omega
^1(M)\times C^\infty (M,\R )\cong \Gamma (T^\ast M\times \R)$ and
$X_0=(-E,0)\in \frak X(M)\times C^\infty (M,\R )\cong \Gamma
(TM\times \R)$. Other relations between generalized Lie bialgebroids
and Jacobi structures were discussed in \cite{IM}.

A generalized Lie bialgebroid $((A,\phi _0),(A^\ast ,X_0))$ is a
generalized Lie bialgebra if the base space $M$ is a single point or,
in other words, if $A$ is a real Lie algebra $\frak g$ of finite
dimension. In \cite{IM}, we started the study of generalized Lie
bialgebras. In particular, we proved that examples of them can be
obtained from algebraic Jacobi structures on a Lie algebra. On the
other hand, we remark that a generalized Lie bialgebra $((\frak g
,\phi _0),(\frak g ^\ast ,X_0))$ such that the 1-cocycles $\phi _0$
and $X_0$ are zero is just a Lie bialgebra in the sense of Drinfeld
\cite{D}. We also recall that there is a one-to-one correspondence
between Lie bialgebras and connected simply connected Poisson Lie
groups (see \cite{D,KM,L,LW,V}).

The aim of this paper is to continue the study of generalized Lie
bialgebras and, more precisely, to discuss some relations between them
and certain types of Jacobi structures on Lie groups. The paper is
organized as follows. In Section 2, we recall several definitions and
results about Jacobi manifolds and generalized Lie bialgebroids which
will be used in the sequel. In Section 3, we show that generalized
Lie bialgebras are closely related with Jacobi structures. In fact,
we prove that the first ones can be considered as the
infinitesimal invariants of Lie groups endowed with special Jacobi
structures (see Theorems \ref{subida} and \ref{bajada}). In Section
4, we propose a method to obtain generalized Lie bialgebras (it is a
generalization of the well-known Yang-Baxter equation method for Lie
bialgebras). As a consequence, we deduce that generalized Lie
bialgebras can be obtained from algebraic Jacobi structures on a Lie
algebra. These results allow us to give, in Section 5, several
examples. In particular, using convenient algebraic contact (locally
conformal symplectic) structures, we obtain interesting examples of
generalized Lie bialgebras. In Section 6, we describe the structure
of a generalized Lie bialgebra $((\frak g ,\phi _0),(\frak g ^\ast
,X_0))$ such that $\frak g$ is a compact Lie algebra and $\phi _0\neq
0$ or $X_0\neq 0$. Finally, the paper closes with two Appendixes. In
the first one, we discuss some relations between algebraic Jacobi
structures and contact (respectively, locally conformal symplectic)
Lie algebras. In the second one, we give a simple proof of the
following assertion: if $\frak h$ is a compact contact Lie algebra of
dimension $\geq 3$ then $\frak h$ is isomorphic to $\frak s\frak u
(2)$. Moreover, we describe all the algebraic contact structures on
$\frak s\frak u (2)$. We use these results in Section 6. We remark
that the aforementioned assertion may be deduced as a corollary of a
more general theorem (about orthogonal contact Lie algebras) which
was proved by Diatta in \cite{Di} (for more details, see Appendix B).

{\bf Notation:} Throughout this paper, we will use the following
notation. If $M$ is a differentiable manifold of dimension $n$, we will denote
by $C^\infty (M,\R)$ the algebra of $C^\infty$ real-valued functions
on $M$, by $\Omega ^k(M)$ the space of $k$-forms on $M$, by $\frak X
(M)$ the Lie algebra of vector fields, by $\delta$ the usual
differential on $\Omega ^\ast (M)=\oplus _k \Omega ^k (M)$ and by
$[\, ,\, ]$ the Schouten-Nijenhuis bracket (\cite{BV,V}). On the
other hand, if $G$ is a Lie group with Lie algebra $\frak g$, we will
denote by $\frak e$ the identity element of $G$, by $L_g:G\to G$
(respectively, $R_g:G\to G$) the left (respectively, right)
translation by $g\in G$, by $Ad: G\times \wedge ^k\frak g\to \wedge
^k\frak g$ the adjoint action of $G$ on $\wedge ^k\frak g$ and by
$ad:\frak g \times \wedge  ^k\frak g\to \wedge ^k\frak g$ the
adjoint representation of $\frak g$ on $\wedge ^k\frak g$, that is,
$ad=T\,Ad$. Moreover, if $s\in \wedge ^k\frak g$ then $\bar{s}$
(respectively, $\tilde{s}$) is the left (respectively, right)
invariant $k$-vector on $G$ defined by $\bar{s}(g)=(L_g)_\ast (s)$
(respectively, $\tilde{s}(g)=(R_g)_\ast (s)$), for all $g\in G$, and
if $\hat{P}$ is a k-vector on $G$ then $\hat{P}_r:G\to
\wedge ^k\frak g$ is the map given by
\begin{equation}\label{3.4'}
\hat{P}_r(g)=(R_{g^{-1}})_\ast (\hat{P}(g)),\, \mbox{ for all }g\in G.
\end{equation}
\section{Generalized$\kern-1.25pt$ Lie$\kern-1.25pt$
bialgebroids$\kern-1.25pt$ and$\kern-1.25pt$ Jacobi structures}  
\setcounter{equation}{0}
\subsection{Jacobi structures and Lie algebroids}
A {\em Jacobi structure} on $M$ is a pair $(\Lambda ,E)$, where $\Lambda$ 
is a 2-vector and $E$ is a vector field on $M$ satisfying the following 
properties:
\begin{equation}\label{ecuaciones}
	      [\Lambda ,\Lambda ]=2E\wedge \Lambda ,\hspace{1cm}
	       [E,\Lambda ]=0.
\end{equation}
The manifold $M$ endowed with a Jacobi structure is called a {\em Jacobi 
manifold}. A bracket of functions (the {\em Jacobi bracket}) is defined by 
\begin{equation}\label{Jacobibracket}
\{ f,g\} =\Lambda (\delta f,\delta g)+fE(g)-gE(f),
\end{equation}        
for all $f,g\in C^\infty (M,\R )$. In fact, the space $C^\infty
(M,\R)$ endowed with the Jacobi bracket is a 
{\em local Lie algebra} in the sense of Kirillov (see \cite{K}).
Conversely, a structure of local Lie algebra on $C^\infty (M,\R)$ defines a
Jacobi structure on $M$ (see \cite{GL,K}). If the vector field $E$ identically
vanishes then $(M,\Lambda )$ is a {\em Poisson manifold}. Jacobi and Poisson 
manifolds were introduced by Lichnerowicz (\cite{Li1,Li2}) (see also
\cite{BV,DLM,K,LM,V,We}). 

A {\em Lie algebroid} $A$ over a manifold $M$ is a vector bundle
$A$ over $M$ together with a Lie bracket $\lcf\, ,\,\rcf$ on the space
$\Gamma (A)$ of the global cross sections of $A\to M$ and a bundle
map $\rho :A \to TM$, called the {\em anchor map}, such that if we
also denote by $\rho :\Gamma (A) \to \frak X (M)$ the homomorphism of
$C^\infty (M,\R )$-modules induced by the anchor map then:
\begin{enumerate}
\item[(i)] $\rho :(\Gamma (A),\lcf \, ,\, \rcf )\to (\frak X
(M),[\, ,\, ])$ is a Lie algebra homomorphism and 
\item[(ii)] for all $f\in C^\infty (M,\R)$ and for all $X ,Y \in
\Gamma (A)$, one has  
	 $$\lcf X ,fY \rcf=f\lcf X,Y\rcf +(\rho  (X )(f))Y.$$
\end{enumerate}
The triple $(A,\lcf \, ,\, \rcf,\rho )$ is called a {\em Lie algebroid
over} $M$ (see \cite{Mk,Pr}).

A real Lie algebra of finite dimension is a Lie algebroid over a
point. Another trivial example of a Lie algebroid is the triple
$(TM,[\, ,\, ],Id)$, where $M$ is a differentiable manifold and
$Id:TM\to TM$ is the identity map.

If $A$ is a Lie algebroid, the Lie bracket on $\Gamma (A)$ can
be extended to the so-called {\em Schouten bracket} $\lcf \, ,\,
\rcf$ on the space $\Gamma (\wedge ^\ast A)=\oplus_k\Gamma (\wedge
^kA)$ of multi-sections of $A$ in such a way that $(\oplus_k\Gamma
(\wedge ^kA),\wedge ,\lcf \, ,\, \rcf )$ is a graded Lie algebra.

On the other hand, imitating the usual differential on the space
$\Omega ^\ast (M)$, we define the {\em differential of the Lie
algebroid} $A$, $d: \Gamma (\wedge ^kA^\ast ) \to \Gamma (\wedge
^{k+1}A^\ast )$, as follows. For $\omega\in \Gamma (\wedge ^kA^\ast
)$ and $X_0,\ldots ,X_k\in \Gamma (A)$,
$$
\begin{array}{ccl}
d\omega (X_0,\ldots,X_k)&=&\displaystyle
\sum_{i=0}^k(-1)^i\rho(X_i)(\omega (X_0,\ldots ,\hat{X}_i,\ldots
,X_k))\\ 
&+&\displaystyle\sum_{i<j}(-1)^{i+j}\omega (\lcf X_i,X_j\rcf
,X_0,\ldots ,\hat{X}_i,\ldots ,\hat{X}_j,\ldots ,X_k).
\end{array}
$$
Moreover, since $d^2=0$, we have the corresponding cohomology spaces.
This cohomology is the {\em Lie algebroid cohomology with
trivial coefficients} (see \cite{Mk}).

Using the above definitions, it follows that a 1-cochain $\phi \in
\Gamma (A^\ast)$ is a 1-cocycle if and only if 
$\phi \lcf X ,Y \rcf =\rho (X )(\phi (Y))-\rho (Y)(\phi
(X ))$, for all $X ,Y \in \Gamma (A)$.

Next, we will consider two examples of Lie algebroids.

{\bf 1.-} {\em The Lie algebroid} $(TM\times \R ,\makebox{{\bf [}\,
,\, {\bf ]}}, \pi)$

If $M$ is a differentiable manifold, then the triple
$(TM\times \R ,\makebox{{\bf [}\, ,\, {\bf ]}}, 
\pi)$ is a Lie algebroid over $M$, where $\pi :TM\times \R \to TM$
is the  canonical projection over the first factor and $\makebox{{\bf
[}\, ,\, {\bf ]}}$ is the bracket given by (see \cite{IM0,IM,NVQ,V2})
\begin{equation}\label{corchetedeprimerorden}
\makebox{{\bf [}} (X,f),(Y,g)\makebox{{\bf ]}}=([X,Y],X(g)-Y(f)), 
\end{equation}
for $(X,f),(Y,g )\in \frak X (M)\times C^\infty (M,\R )\cong \Gamma
(TM\times \R)$. 

{\bf 2.-} {\em The Lie algebroid $(T^\ast M \kern-1pt \times \R ,\lcf ,\rcf
_{(\Lambda ,E)},\kern-2pt \widetilde{\#}_{(\Lambda ,E)})$ associated with
a Jacobi manifold \kern-1pt $(\kern-2pt M,\Lambda ,E)$}

If $A\to M$ is a vector bundle over $M$ and $P \in
\Gamma (\wedge ^2 A)$ is a 2-section of $A$, we will denote by
$\#_P:\Gamma (A^\ast)\to \Gamma (A)$ the homomorphism of $C^\infty
(M,\R)$-modules given by $\beta (\#_P(\alpha ))=P (\alpha ,\beta )$,
for $\alpha, \beta \in \Gamma (A^\ast )$. 
We will also denote by $\#_P:A^\ast\to A$ the corresponding bundle
map. Then, a Jacobi manifold $(M,\Lambda ,E)$ has an associated Lie algebroid
$(T^\ast M \times \R ,\lcf ,\rcf _{(\Lambda ,E)},$ $\widetilde{\#}_{(\Lambda
,E)})$, where $\lcf ,\rcf _{(\Lambda ,E)}$ and
$\widetilde{\#}_{(\Lambda,E)}$ are defined by  
\begin{equation}\label{ecjacobi}
\begin{array}{cll}
\kern-17.6pt\lcf (\alpha ,f),(\beta ,g)\rcf _{(\Lambda
,E)}&\kern-9pt=&\kern-9pt({\cal L}_{\#_{\Lambda}(\alpha
)}\beta\kern-2pt-\kern-2pt{\cal L}_{\#_{\Lambda}(\beta
)}\alpha\kern-2pt -\kern-2pt \delta (\Lambda
(\alpha ,\beta ))\kern-2pt+\kern-2ptf{\cal
L}_{E}\beta\kern-2pt-\kern-2ptg {\cal L}_{E}\alpha 
\kern-2pt-\kern-2pti(E)(\alpha \wedge \beta),\\ 
&\kern-9pt &\kern-9pt \Lambda (\beta ,\alpha
)\kern-2pt+\kern-2pt\#_{\Lambda}(\alpha
)(g)\kern-2pt-\kern-2pt\#_{\Lambda}(\beta
)(f)\kern-2pt+\kern-2ptfE(g)\kern-2pt-\kern-2ptg E(f)),\\ 
& & \\
\kern-17.6pt\widetilde{\#}_{(\Lambda ,E)}(\alpha
,f)&\kern-9pt=&\#_{\Lambda}(\alpha )+fE, 
\end{array}
\end{equation}
for $(\alpha ,f),(\beta ,g)\in \Omega ^1(M)\times
C^\infty (M,\R )$, ${\cal L}$ being the Lie derivative operator (see
\cite{KS}). 

In the particular case when $(M,\Lambda )$ is a Poisson manifold we
recover, by projection, the Lie algebroid $(T^\ast M,\lcf \, ,\, \rcf
_{\Lambda}, \#_{\Lambda})$, where $\lcf \, ,\, \rcf_{\Lambda}$ is the
bracket of 1-forms defined by $\lcf \alpha ,\beta \rcf _{\Lambda}={\cal
L}_{\#_{\Lambda}(\alpha )}\beta-{\cal L}_{\#_{\Lambda}(\beta )}\alpha
-\delta (\Lambda (\alpha ,\beta )),$ for $\alpha ,\beta \in \Omega ^1(M)$
(see \cite{BV,CDW,F,V}).
\subsection{Generalized Lie bialgebroids}
In this Section, we will recall the definition of a generalized Lie
bialgebroid. First, we will exhibit some results about the
differential calculus on Lie algebroids in the presence of a
1-cocycle (for more details, see \cite{IM}).

If $(A,\lcf \, ,\, \rcf ,\rho )$ is a Lie algebroid over $M$ and, in
addition, we have a 1-cocycle $\phi _0\in \Gamma (A ^\ast)$ then the
usual representation of the Lie algebra $\Gamma (A)$ on the space
$C^\infty (M,\R )$ can be modified and a new representation is
obtained. This representation is given by $\rho _{\phi _0}(X)(f)=\rho
(X)(f)+\phi _0 (X)f$, for $X\in \Gamma (A)$ and $f\in C^\infty (M,\R )$.
The resultant cohomology operator $d_{\phi _0}$ is called
the $\phi _0$-differential of $A$ and its expression, in terms of the
differential $d$ of $A$, is $d_{\phi _0}\omega =d\omega +\phi
_0\wedge \omega$, for $\omega \in \Gamma (\wedge ^k A^\ast )$. The
$\phi _0$-differential of $A$ allows us to define, in a natural way,
the $\phi _0$-Lie derivative by a section $X\in \Gamma (A)$, $({\cal
L}_{\phi _0})_X:\Gamma (\wedge ^kA^\ast)\to \Gamma (\wedge
^kA^\ast)$, as the commutator of $d_{\phi _0}$ and the 
contraction by $X$, that is, $({\cal L}_{\phi _0})_X  =d_{\phi
_0}\circ i(X)+i(X)\circ d_{\phi _0}$ (for the general definition of
the differential and the Lie derivative associated with a
representation of a Lie algebroid on a vector bundle, see \cite{Mk}). 

On the other hand, imitating the definition of the Schouten bracket
of two multilinear first-order differential operators on the space of
$C^\infty$ real-valued functions on a manifold $N$ (see \cite{BV}),
we introduced the $\phi _0$-Schouten bracket of a $k$-section $P$ and
a $k'$-section $P'$ as the $k+k'-1$-section given by 
\begin{equation}\label{Schouext}
\lcf P, P' \rcf_{\phi _0}=\lcf P, P' \rcf + (-1)^{k+1}(k-1)P\wedge
(i(\phi _0)P')-(k'-1)(i(\phi_0)P)\wedge P',
\end{equation}
where $\lcf \, ,\, \rcf$ is the usual Schouten bracket of $A$. The
$\phi _0$-Schouten bracket satisfies the following properties. For 
$f\in C^\infty (M,\R )$, $X ,Y \in \Gamma (A)$, $P \in \Gamma
(\wedge ^k A)$, $P' \in \Gamma (\wedge ^{k'} A)$ and $P'' \in \Gamma
(\wedge ^{k''} A)$,
$$
\lcf X ,f\rcf _{\phi _0} = \rho _{\phi _0}(X )(f),\quad \lcf X ,Y
\rcf _{\phi _0}=\lcf X ,Y \rcf,\quad \lcf P, P'\rcf _{\phi
_0}=(-1)^{kk'}\lcf P', P\rcf _{\phi _0}, 
$$
$$
\lcf P, P'\wedge P''\rcf _{\phi _0}=\lcf P, P'\rcf _{\phi _0}\wedge P''
+(-1)^{k'(k+1)} P'\wedge \lcf P, P''\rcf _{\phi _0}-(i(\phi
_0)P)\wedge  P'\wedge P'',
$$
$$
   (-1)^{kk''}\lcf \lcf P,P' \rcf _{\phi _0},P''\rcf _{\phi _0}+
   (-1)^{k'k''}\lcf \lcf P'',P \rcf _{\phi _0},P'\rcf _{\phi _0}+
   (-1)^{kk'}\lcf \lcf P',P'' \rcf _{\phi _0},P\rcf _{\phi _0}=0.
$$   
Using the $\phi _0$-Schouten bracket, we can
define the $\phi _0$-Lie derivative of $P\in \Gamma (\wedge ^k A)$ by
$X \in \Gamma (A)$ as $({\cal L}_{\phi _0})_X (P)=\lcf X ,P\rcf _{\phi _0}.$

Now, suppose that $(A,\lcf \, ,\, \rcf ,\rho)$ is a Lie algebroid and
that $\phi _0\in \Gamma (A^\ast )$ is a 1-cocycle. Assume also that
the dual bundle $A^\ast$ admits a Lie algebroid structure $(\lcf \,
,\, \rcf _\ast ,\rho _\ast )$ and that $X_0\in \Gamma
(A)$ is a 1-cocycle. The pair $((A,\phi _0),(A^\ast ,X_0))$ is a {\em
generalized Lie bialgebroid} if  
\begin{equation}\label{condcomp}
 d_\ast{}_{X_0}\lcf X ,Y \rcf = \lcf X ,d_\ast{}_{X_0}Y \rcf _{\phi
_0}-\lcf Y ,d_\ast{}_{X_0}X \rcf _{\phi _0},\quad
 ({\cal L}_\ast {}_{X_0})_{\phi _0}P+({\cal L}_{\phi _0})_{X_0}P= 0,
\end{equation}
for $X ,Y \in \Gamma (A)$ and $P\in \Gamma (\wedge ^kA)$, where
$d_\ast{}_{X_0}$ (respectively, ${\cal L}_\ast{}_{X_0}$) is the
$X_0$-differential (respectively, the $X_0$-Lie derivative) of 
$A^\ast$. Note that the second equality in (\ref{condcomp}) holds if
and only if $\phi _0(X_0)=0$, $\rho (X_0)=-\rho _\ast (\phi _0)$ and
$({\cal L}_\ast {}_{X_0})_{\phi _0}X+\lcf X_0,X\rcf = 0$,
for $X\in \Gamma (A)$ (for more details, see \cite{IM}). Moreover, in
the particular case when $\phi _0=0$ and $X_0=0$, (\ref{condcomp}) is
equivalent to the condition $ d_\ast{}\lcf X ,Y \rcf = \lcf X ,d_\ast
Y \rcf -\lcf Y ,d_\ast X \rcf .$ 
Thus, the pair $((A,0),(A^\ast ,0))$ is a generalized Lie bialgebroid
if and only if the pair $(A,A^\ast)$ is a Lie bialgebroid (see
\cite{K-S,MX}). 

On the other hand, if $(M,\Lambda ,E)$ is a Jacobi manifold, then we
proved in \cite{IM} that the pair $\Big ((TM\times \R ,\phi
_0)$,$(T^\ast M \times \R ,X_0)\Big )$ is 
a generalized Lie bialgebroid, where $\phi _0$ and $X_0$ are the
1-cocycles on $TM\times \R$ and $T^\ast M\times \R$ given by
$$\phi _0=(0,1)\in \Omega ^1(M)\times C^\infty (M,\R )\cong \Gamma
(T^\ast M\times \R),$$  
$$X_0=(-E,0)\in \frak X (M)\times C^\infty (M,\R )\cong \Gamma
(TM\times \R).$$  

\section{Generalized Lie bialgebras and Jacobi structures on
connected Lie groups}
\setcounter{equation}{0}
In this Section, we will deal with generalized Lie bialgebroids over a
point. 
\begin{definition}\cite{IM}
A generalized Lie bialgebra is a generalized Lie bialgebroid over a
point, that is, a pair $((\frak g ,\phi _0),(\frak g ^\ast ,X_0))$,
where $(\frak g ,[\, ,\, ]^\frak g )$ is a real Lie algebra of finite
dimension such that the dual space $\frak g ^\ast$ is also a Lie
algebra with Lie bracket $[\, ,\, ]^{\frak g ^\ast}$, $X_0 \in \frak g$ and
$\phi _0\in \frak g^\ast$ are 1-cocycles on $\frak g ^\ast$ and
$\frak g$, respectively, and 
\begin{equation}\label{condalg1}
d_\ast {}_{X_0}[X,Y]^\frak g =[X,d_\ast{}_{X_0}Y]^\frak g _{\phi
_0}-[Y,d_\ast{}_{X_0}X]^\frak g _{\phi _0},
\end{equation}
\begin{equation}\label{condalg2}
\phi _0 (X_0)=0,
\end{equation}
\begin{equation}\label{condalg3}
i(\phi _0)(d_\ast X)+[X_0,X]^\frak g =0,
\end{equation}
for all $X,Y\in \frak g$, $d_\ast$ being the algebraic differential
on $(\frak g ^\ast ,[\, ,\, ]^{\frak g ^\ast})$.
\end{definition}
\begin{remark}
{\rm
In the particular case when $\phi _0=0$ and $X_0=0$, we recover the
concept of a {\em Lie bialgebra}, that is, a dual pair $(\frak g,\frak g
^\ast)$ of Lie algebras such that $d_\ast{}[ X ,Y ]^\frak g = [ X
,d_\ast Y ]^\frak g -[ Y ,d_\ast X]^\frak g ,$ for $X,Y\in \frak g$
(see \cite{D}). 
}
\end{remark}
We know that there exists a one-to-one correspondence between Lie
bialgebras and connected simply connected Poisson Lie groups (see
\cite{D,L,LW,V}). So, we will study a connected Lie group $G$ with Lie
algebra $\frak g$ such that the pair $((\frak g ,\phi _0),(\frak g
^\ast ,X_0))$ is a generalized Lie bialgebra. 

We will use the following well-known results about cocycles on
Lie groups and on their Lie algebras.
\begin{lemma}\label{Lu}\cite{L}
Let $G$ be a connected Lie group with Lie algebra $\frak g$. Let
$\Phi :G\times V\to V$ be a representation of $G$ on a vector space
$V$. Let $T\Phi:\frak g \times V\to V$ be the induced representation
of $\frak g$ on $V$. 
\begin{itemize}
\item[{\it i)}] If the map $\phi :G\to V$ is a 1-cocycle on $G$
relative to $\Phi$, i.e., if for $h,g\in G$
       $$\phi (hg)=\phi (h)+\Phi (h, \phi (g)),$$
then $\epsilon=:(\delta \phi)(\frak e):\frak g \to V$, the derivative of
$\phi$ at $\frak e$, is a 1-cocycle on $\frak g$ relative to $T\Phi$, i.e.,
for $X,Y\in \frak g$
  $$T\Phi (X,\epsilon (Y))-T\Phi (Y, \epsilon (X))=\epsilon
			   ([X,Y]^\frak g).$$
Moreover, $\delta \phi =0$ implies that $\phi =0$.
\item[{\it ii)}] When $G$ is simply connected, any 1-cocycle
$\epsilon$ on $\frak g$ relative to $T\Phi$ can be integrated to give
a unique 1-cocycle $\phi$ on $G$ relative to $\Phi$ such that
$(\delta \phi)(\frak e) =\epsilon$. 
\item[{\it iii)}] When $\frak g$ is semisimple, every 1-cocycle
$\epsilon :\frak g \to V$ on $\frak g$ is a coboundary, that is,
$\epsilon (X)=T\Phi (X, v_0)$, for some $v_0\in V$. 
\end{itemize}
\end{lemma}
Next, we will introduce the notion of a ($\sigma ,c$)-multiplicative
$k$-vector on a connected Lie group $G$, where $\sigma :G\to \R$ is a
multiplicative function and $c\in\R$. This notion will play an
important role (in Section \ref{Correspondencia}) in the description
of the Jacobi structure on a connected simply connected
Lie group whose Lie algebra is $\frak g$ and such that the pair
$((\frak g,\phi _0),(\frak g^\ast ,X_0))$ is a generalized Lie
bialgebra. We recall that a C$^\infty$ real-valued function $\sigma
:G\to \R$ is multiplicative if it is a Lie group
homomorphism. 
\subsection{($\sigma ,c$)-multiplicative multivectors on a Lie group}
We will denote by $e:\R\to \R$ the real exponential. Then,
\begin{proposition}\label{equivalencias}
Let $G$ be a connected Lie group, $\sigma :G\to \R$ a
multiplicative function and $c\in\R$. If $\hat{P}$ is a k-vector on
$G$, the following properties are equivalent:
\begin{itemize}
\item[{\it i)}] $\hat{P}_r (hg)=\hat{P}_r(h)+e^{-(k-c)\sigma (h)}Ad
_h(\hat{P}_r(g))$. 
\item[{\it ii)}] $\hat{P}(hg)=(R_g)_\ast (\hat{P}(h))+e^{-(k-c)\sigma
(h)}(L_h)_\ast (\hat{P}(g))$.  
\item[{\it iii)}] $e^{(k-c)\sigma (hg)}\hat{P}(hg)=e^{(k-c)\sigma
(g)}(R_g)_\ast (e^{(k-c)\sigma (h)}\hat{P}(h))+(L_h)_\ast
(e^{(k-c)\sigma (g)}\hat{P}(g))$. 
\item[{\it iv)}] $\hat{P}(\frak e)=0$ and $e^{(k-c)\sigma}{\cal
L}_{\bar{X}}\hat{P}$ is left invariant whenever $\bar{X}$ is a left
invariant vector field on $G$.
\item[{\it v)}] $\hat{P}(\frak e)=0$ and $e^{-(k-c)\sigma}{\cal
L}_{\tilde{X}}(e^{(k-c)\sigma}\hat{P}) $ is right invariant whenever
$\tilde{X}$ is a right invariant vector field on $G$. 
\end{itemize}
\end{proposition}
\prueba The result follows using (\ref{3.4'}), the fact that $\sigma$ is a
multiplicative function and proceeding as in the proof of Proposition
10.5 in \cite{V}.\QED 

Now, we introduce the definition of a $(\sigma ,c)$-multiplicative
$k$-vector on $G$.
\begin{definition}\label{sigma-mult} 
Let $G$ be a connected Lie group, $\sigma :G\to \R$ a
multiplicative function and $c\in\R$. A k-vector $\hat{P}$ on $G$ is
said to be ($\sigma ,c$)-multiplicative if $\hat{P}$ satisfies any of
the properties in Proposition \ref{equivalencias}. In particular, if
$c=1$, we will say that the k-vector is $\sigma$-multiplicative.
\end{definition}
It is clear that if $\hat{P}$ is a $(\sigma ,c$)-multiplicative k-vector
and $\sigma$ identically vanishes, then $\hat{P}$ is multiplicative
(see \cite{L,V}). 

Let $G$ be a connected Lie group with Lie algebra $\frak g$,
$\sigma :G\to \R$ a multiplicative function and $c\in\R$. We can
introduce the representation $Ad_{(\sigma ,c)}:G\times \wedge ^k\frak g\to 
\wedge ^k \frak g$ of $G$ on $\wedge ^k\frak g$ defined by
\begin{equation}\label{3.4iv}
(Ad_{(\sigma ,c)})_g(s)= e^{-(k-c)\sigma (g)}Ad_gs,
\end{equation}
for $g\in G$ and $s\in \wedge ^k\frak g$. If $\phi _0=(\delta \sigma
)(\frak e)$ then we will denote by $ad_{(\phi _0 ,c)}$ the corresponding
representation of $\frak g$ on $\wedge ^k\frak g$, that is, $ad_{(\phi
_0,c)}=T\, Ad_{(\sigma ,c)}:\frak g \times \wedge ^k\frak g\to \wedge ^k\frak
g$. From (\ref{3.4iv}), it follows that
\begin{equation}\label{3.8'}
ad_{(\phi _0,c)}(X)(s)=[X,s]^\frak g -(k-c)\phi
_0(X)s=ad(X)(s)-(k-c)\phi _0(X)s
\end{equation}
for $X\in \frak g$ and $s\in \wedge ^k\frak g$, where $[\, ,\,
]^\frak g$ is the Schouten bracket of the Lie algebroid $\frak g\to
\{ \mbox{ a single point }\}$. It is clear that $\phi _0\in \frak g
^\ast$ is a 1-cocycle with respect to the trivial representation of
$\frak g$ on $\R$ and that if $c=1$ then (see (\ref{Schouext})) 
\begin{equation}\label{3.8''}
ad_{(\phi _0,1)}(X)(s)=[X,s]^\frak g_{\phi _0}.
\end{equation}
\begin{remark}\label{3.6'}
{\rm
Note that if $\hat{P}$ is a k-vector on $G$ then, from (\ref{3.4iv})
and Proposition \ref{equivalencias}, we obtain that $\hat{P}$ is $(\sigma
,c)$-multiplicative if and only if $\hat{P}_r:G\to \wedge ^k\frak g$
is a 1-cocycle with respect to the representation $Ad_{(\sigma
,c)}:G\times \wedge ^k\frak g\to \wedge ^k\frak g$.
}
\end{remark}
Now, suppose that $\hat{P}$ is a k-vector on $G$ such that
$\hat{P}(\frak e)=0$. Then, one can define the intrinsic derivative of
$\hat{P}$ at $\frak e$ as the linear map $\delta _\frak e\hat{P}:\frak g\to
\wedge ^k\frak g$ given by (see \cite{L,V})
\begin{equation}\label{3.4v}
(\delta _\frak e\hat{P})(X)=(\delta \hat{P}_r)(\frak e)(X)=({\cal
L}_{\hat{X}}\hat{P})(\frak e),
\end{equation}
for $X\in \frak g$, $\hat{X}$ being an arbitrary vector field on $G$
satisfying $\hat{X}(\frak e)=X$. Using (\ref{3.4v}), Lemma
\ref{Lu} and Remark \ref{3.6'}, we deduce
\begin{proposition}\label{nuestroLu}
Let $G$ be a connected Lie group, $\sigma :G\to \R$ a
multiplicative function and $c\in\R$. Suppose that $\phi _0=(\delta
\sigma )(\frak e)$. 
\begin{itemize}
\item[{\it i)}] If $\hat{P}$ is a $(\sigma ,c)$-multiplicative k-vector
then its intrinsic derivative $\delta _\frak e\hat{P}:\frak g\to
\wedge ^k\frak g$ is a 1-cocycle with respect to the representation
$ad_{(\phi _0,c)}:\frak g\times \wedge ^k\frak g\to \wedge ^k\frak g$.
\item[{\it ii)}] If $G$ is simply connected and $\epsilon :\frak
g\to \wedge ^k\frak g$ is a 1-cocycle with respect to the representation
$ad_{(\phi _0,c)}:\frak g\times \wedge ^k\frak g\to \wedge ^k\frak g$
then there exists a unique $(\sigma ,c)$-multiplicative k-vector $\hat{P}$
such that its intrinsic derivative at $\frak e$, $\delta _\frak
e\hat{P}$, is just $\epsilon$. 
\end{itemize}
\end{proposition}
\begin{remark}\label{nota3.8'}
{\rm
Let $G$ be a connected Lie group, $\sigma :G\to\R$ a
multiplicative function and $c\in\R$. If $\hat{P}$ is a $(\sigma
,c)$-multiplicative $k$-vector then, from Proposition
\ref{equivalencias}, it follows that
$$(\delta \hat{P}_r)(h)((L_h)_\ast (X))=e^{-(k-c)\sigma
(h)}Ad_h((\delta _\frak e\hat{P})(X)),$$ 
for $h\in G$ and $X\in \frak g$. Thus, $\delta \hat{P}_r=0$ if and
only if the intrinsic derivative of $\hat{P}$ at $\frak e$ is zero.
Therefore, $\hat{P}=0$ if and only if the intrinsic derivative
of $\hat{P}$ at $\frak e$ is zero (see Lemma \ref{Lu} and Remark \ref{3.6'}).
}
\end{remark}
\begin{example}\label{cobordes}
{\rm
Let $G$ be a connected Lie group with Lie algebra $\frak g$, $\sigma
:G\to \R$ a multiplicative function and $c\in\R$. Suppose that $s\in\wedge 
^k\frak g$. Then, we consider the k-vector $\hat{s}$ on $G$ defined by
$$\hat{s}(g)=e^{-(k-c)\sigma (g)}\bar{s}(g)-\tilde{s}(g),\mbox{ for
all }g\in G.$$
A direct computation shows that $\hat{s}$ is a
$(\sigma ,c)$-multiplicative k-vector on $G$. Moreover, the intrinsic
derivative of $\hat{s}$ at $\frak e$ is given by
 $$(\delta _\frak e \hat{s})(X)=[X,s]^\frak g -(k-c)\phi _0(X)s=ad_{(\phi
_0,c )}(X)(s),$$ 
for $X\in \frak g$, where $\phi _0=(\delta\sigma )(\frak e)$. Note that, in
this case, $\delta _\frak e \hat{s}$ is a 1-coboundary with respect to the
representation $ad_{(\phi _0,c)}:\frak g\times \wedge ^k\frak g \to
\wedge ^k\frak g$. Moreover, using Remark \ref{nota3.8'}, we deduce
that $s$ is $ad_{(\phi _0,c)}$-invariant if and only if
$\tilde{s}=e^{-(k-c)\sigma }\bar{s}$. 
}
\end{example}
\subsection{Generalized Lie bialgebras and Jacobi structures on
connected Lie groups}\label{Correspondencia}
We will prove that if $G$ is a connected simply connected Lie group
with Lie algebra $\frak g$ and the pair $((\frak g,\phi _0),(\frak g
^\ast ,X_0))$ is a generalized Lie bialgebra then $G$ admits an
special Jacobi structure.
\begin{theorem}\label{subida}
Let $((\frak g,\phi _0),(\frak g ^\ast ,X_0))$ be a generalized Lie
bialgebra and $G$ a connected simply connected Lie group with Lie
algebra $\frak g$. Then, there exists a unique multiplicative function
$\sigma :G\to \R$ and a unique $\sigma$-multiplicative 2-vector
$\Lambda$ on $G$ such that $(\delta \sigma )(\frak e)=\phi _0$ and
the intrinsic derivative of $\Lambda$ at $\frak e$ is $-d_{\ast
X_0}$. Moreover, the following relation holds
\begin{equation}\label{compatibilidad}
\# _\Lambda (\delta \sigma )=\tilde{X}_0-e^{-\sigma}\bar{X}_0,
\end{equation}
and the pair $(\Lambda ,E)$ is a Jacobi structure on
$G$, where $E=-\tilde{X}_0$.
\end{theorem}
\prueba Since $G$ is connected and simply connected then, using Lemma
\ref{Lu} and the fact that $\phi _0$ is a 1-cocycle with respect to
the trivial representation of $G$ on $\R$, we deduce that there exists
a unique multiplicative function $\sigma :G\to\R$ satisfying $(\delta \sigma
)(\frak e)=\phi _0$. Now, take $\epsilon :\frak g\to \wedge ^2\frak
g$ given by $\epsilon (X)=-d_\ast{}_{X_0}X$, for $X\in \frak g$. From
(\ref{condalg1}), it follows that $\epsilon$ is a 1-cocycle of $\frak
g$ with respect to the representation $ad_{(\phi _0 ,1)}$. Thus,
using Proposition \ref{nuestroLu}, we have that there exists a unique
$\sigma$-multiplicative 2-vector $\Lambda$ on $G$ such that its intrinsic 
derivative at $\frak e$ is just $-d_\ast{}_{X_0}$, that is,
\begin{equation}\label{3.10'}
\delta _\frak e\Lambda =-d_\ast{}_{X_0}.
\end{equation}
Next, we will see that (\ref{compatibilidad}) holds. Let $\bar{X}\in
\frak X (G)$ be a left invariant vector field and $X=\bar{X}(\frak e)$.
Then, ${\cal L}_{\bar{X}}(\delta \sigma )=\delta ({\cal
L}_{\bar{X}}\sigma )=0$ and
$$
[\bar{X},\# _{\Lambda}(\delta \sigma
)+e^{-\sigma}\bar{X}_0]=i(\delta \sigma)({\cal L}_{\bar{X}}\Lambda )
-e^{-\sigma}\phi _0(X)\bar{X}_0+e^{-\sigma}\overline{[X,X_0]^\frak g}.
$$
The 2-vector $e^\sigma {\cal L}_{\bar{X}}\Lambda$ is left invariant
(see Proposition \ref{equivalencias}). Therefore, if $g\in G$ and
$\alpha _g\in T^\ast _gG$, we obtain that 

$\begin{array}{lcl}
\alpha _g \Big \{ i(\delta \sigma )({\cal L}_{\bar{X}}\Lambda
)-e^{-\sigma}\phi _0(X)\bar{X}_0+e^{-\sigma}\overline{[X,X_0]^\frak
g}\Big \} (g)&=& \\
&\kern-325pt=&\kern-150pte^{-\sigma (g)}\Big ( ({\cal
L}_{\bar{X}}\Lambda )_{(\frak e)}  
(((L_g)_\ast )^\ast ((\delta\sigma)(g)),((L_g)_\ast )^\ast \alpha _g)\\
&\kern-325pt&\kern-150pt-\phi _0 (X)(((L_g)_\ast )^\ast (\alpha
_g))(X_0) +(((L_g)_\ast )^\ast (\alpha _g))([X,X_0]^\frak g )\Big ), 
\end{array}
$

where $((L_g)_\ast )^\ast :T^\ast _gG\to \frak g ^\ast$ is the
adjoint homomorphism of $(L_g)_\ast :\frak g\to T_gG$.

Note that the 1-form $\delta \sigma$ is left invariant which implies
that $((L_g)_\ast )^\ast (\delta\sigma (g))=\phi _0$. Consequently,
from (\ref{condalg2}), (\ref{condalg3}), (\ref{3.4v}) and
(\ref{3.10'}), we deduce that 
$$\alpha _g \Big \{ i(\delta \sigma )({\cal L}_{\bar{X}}\Lambda
)-e^{-\sigma}\phi _0(X)\bar{X}_0+e^{-\sigma}\overline{[X,X_0]^\frak g
}\Big \} (g)=0.$$
Thus, $\# _\Lambda (\delta \sigma )+e^{-\sigma}\bar{X}_0$ is a
right invariant vector field and, since 
$\Big (\# _\Lambda (\delta \sigma )+e^{-\sigma}\bar{X}_0\Big )(\frak
e)=X_0,$ we conclude that (\ref{compatibilidad}) holds.    

Now, take $E=-\tilde{X}_0$. Since $E$ is a right invariant
vector field and $\Lambda$ is $\sigma$-multiplicative, we have that 
$e^{-\sigma}{\cal L}_E(e^\sigma \Lambda )$ is right invariant (see
Proposition \ref{equivalencias}). Moreover, 
from (\ref{condalg2}), 
   $$e^{-\sigma}{\cal L}_E(e^\sigma \Lambda )=e^{-\sigma}(e^\sigma
  E(\sigma )\Lambda+e^\sigma {\cal L}_E\Lambda )={\cal L}_E\Lambda .$$  
On the other hand, using (\ref{3.4v}), (\ref{3.10'}) and the fact
that $X_0\in\frak g$ is a 1-cocycle (that is, $d_\ast X_0=0$), it
follows that $({\cal L}_E\Lambda )(\frak e)=0$. This implies that ${\cal
L}_E\Lambda =[E,\Lambda ]=0$. 

Finally, we will prove that $[\Lambda ,\Lambda ]-2E\wedge \Lambda
=0$. First, we will show that $[\Lambda ,\Lambda ]-2E\wedge \Lambda$
is $\sigma$-multiplicative. Since $\Lambda (\frak e)=0$, we have that
$([\Lambda ,\Lambda ]-2E\wedge \Lambda )(\frak e)=0$ (see \cite{V}).
Moreover, if $\bar{X}$ is a left invariant vector field then
$[\bar{X},E]=0$ and, using (\ref{compatibilidad}) and the properties
of the Schouten-Nijenhuis bracket, we deduce that
$$
e^{2\sigma}{\cal L}_{\bar{X}}\Big ( [\Lambda ,\Lambda ]-2E\wedge
\Lambda \Big )=2\Big ( e^\sigma [e^\sigma {\cal L}_{\bar{X}}\Lambda
,\Lambda ]+\bar{X}_0\wedge (e^\sigma {\cal L}_{\bar{X}}\Lambda )\Big
) .
$$
On the other hand, from Proposition \ref{equivalencias}, it follows
that $e^\sigma {\cal L}_{\bar{X}}\Lambda$ and $e^\sigma [e^\sigma
{\cal L}_{\bar{X}}\Lambda ,\Lambda ]$ are left invariant
multivectors. Therefore,  $e^{2\sigma}{\cal L}_{\bar{X}}\Big (
[\Lambda ,\Lambda ]-2E\wedge \Lambda \Big )$ is also a left invariant
multivector. Consequently, $[\Lambda ,\Lambda ]-2E\wedge \Lambda$ is
$\sigma$-multiplicative, as we wanted to prove.

Next, we will compute the intrinsic derivative at $\frak e$ of the 3-vector
$[\Lambda ,\Lambda ]-2E\wedge \Lambda$.

If $[\, ,\, ]_\Lambda :\wedge ^2\frak g ^\ast\to\frak g ^\ast$ is the
adjoint map of the intrinsic derivative of $\Lambda$ at $\frak e$, using
(\ref{3.10'}), we obtain that
\begin{equation}\label{3.10''}
[\alpha ,\beta ]_\Lambda =[\alpha ,\beta ]^{\frak g
^\ast}-i(X_0)(\alpha \wedge \beta ),
\end{equation}
for $\alpha ,\beta \in \frak g ^\ast$, where $[\, ,\, ]^{\frak g
^\ast}$ is the Lie bracket on $\frak g ^\ast$. This implies that
\begin{equation}\label{3.10'''}
[\alpha ,\beta ]_\Lambda (X_0)=[\alpha ,\beta ]^{\frak g ^\ast}(X_0)=0.
\end{equation}
Now, from (\ref{3.4v}), (\ref{3.10'}) and since $E=-\tilde{X_0}$, we
have that 
\begin{equation}\label{3.10iv}
\delta _\frak e ([\Lambda ,\Lambda ]-2E\wedge \Lambda )(X)={\cal
L}_{\bar{X}}([\Lambda, \Lambda ]-2E\wedge \Lambda )(\frak e)
=\delta _\frak e[\Lambda ,\Lambda ](X)-2X_0\wedge d_\ast X ,
\end{equation}
for $X\in\frak g$. Thus, using
(\ref{3.10''}), (\ref{3.10'''}) and (\ref{3.10iv}), we conclude that
$$
\Big \{ \delta _\frak e ([\Lambda ,\Lambda ]-2E\wedge \Lambda )(X)\Big \}
(\alpha ,\beta ,\gamma )
=-2\displaystyle\sum_{Cycl.(\alpha ,\beta
,\gamma )}\Big ([\alpha ,[\beta ,\gamma ]^{\frak g ^\ast} ]^{\frak g^\ast}
\Big )(X)=0,
$$
for $\alpha ,\beta ,\gamma \in \frak g ^\ast$, that is, the intrinsic
derivative of $[\Lambda ,\Lambda ]-2E\wedge \Lambda$ at $\frak e$ is null.
Therefore, $[\Lambda ,\Lambda ]=2E\wedge \Lambda$ (see Remark
\ref{nota3.8'}).\QED 
\begin{remark}
{\rm
\begin{itemize}
\item[{\it i)}] Under the same hypotheses as in Theorem \ref{subida},
if $\phi _0=0$ then the multiplicative function $\sigma$ will be
null, the 2-vector $\Lambda$ will be multiplicative and the vector
field $E$ will be bi-invariant (see (\ref{compatibilidad})).
\item[{\it ii)}] Under the same hypotheses as in Theorem
\ref{subida}, if $\phi _0=0$ and $X_0=0$ then $\sigma$ and $E$ will be
null and $(G,\Lambda )$ will be a Poisson Lie group.
\end{itemize}
}
\end{remark}
Now, we discuss a converse of Theorem \ref{subida}.
\begin{theorem}\label{bajada}
Let $(\Lambda ,E)$ be a Jacobi structure on a connected Lie group $G$
and $\sigma :G\to \R$ a multiplicative function such that:
\begin{itemize}
\item[{\it i)}] $\Lambda$ is $\sigma$-multiplicative.
\item[{\it ii)}] $E$ is a right invariant vector field, $E(\frak e)=-X_0$
and $\# _\Lambda (\delta \sigma )=\tilde{X}_0-e^{-\sigma}\bar{X}_0.$
\end{itemize}
If $[\, ,\, ]_\Lambda :\wedge ^2\frak g ^\ast \to \frak g ^\ast$ is
the adjoint map of the intrinsic derivative of $\Lambda$ at $\frak e$ and
$[\, ,\, ]^{\frak g ^\ast}$ is the bracket on $\frak g ^\ast$ given
by  
\begin{equation}\label{3.10v}
[\alpha ,\beta ]^{\frak g ^\ast}=[\alpha ,\beta ]_\Lambda+
i(X_0)(\alpha \wedge \beta ),\mbox{ for }\alpha ,\beta \in \frak g
^\ast ,
\end{equation}
then $(\frak g ^\ast ,[\, ,\, ]^{\frak g ^\ast})$ is a Lie algebra
and the pair $((\frak g ,\phi _0),(\frak g^\ast ,X_0))$ is
a generalized Lie bialgebra, where $\phi _0=(\delta \sigma )(\frak e)$.
\end{theorem}
\prueba Since $\sigma$ is a multiplicative function, we have that
$\phi _0$ is a 1-cocycle of $(\frak g,[\, ,\, ]^{\frak g})$. Now,
suppose that $\alpha ,\beta\in \frak g ^\ast$. We consider 
two $C^\infty$ real-valued functions $f$ and $g$ on $G$ such that
\begin{equation}\label{3.10vi}
f(\frak e)=g(\frak e)=0,\quad (\delta f)(\frak e)=\alpha,\quad
(\delta g)(\frak e)=\beta .
\end{equation}
If $\{ \, ,\, \}$ is the Jacobi bracket associated with the Jacobi
structure $(\Lambda ,E)$ then, from (\ref{Jacobibracket}),
(\ref{3.4v}), (\ref{3.10v}) and (\ref{3.10vi}), we deduce that
\begin{equation}\label{3.10vii}
(\delta \{ f,g\} )(\frak e)=\delta (\Lambda (\delta f, \delta
g))(\frak e)+i(X_0)(\alpha \wedge \beta )
=[\alpha ,\beta ]^{\frak g ^\ast}.
\end{equation}
Using (\ref{3.10vii}) it follows that $(\frak g ^\ast ,[\, ,\,
]^{\frak g^\ast})$ is a Lie algebra. Moreover, from (\ref{3.10v}), we
obtain that 
\begin{equation}\label{3.10viii}
(\delta _\frak e \Lambda)(X)=-d_\ast{}_{X_0}X,\mbox{ for }X\in \frak g .
\end{equation}
Thus, using (\ref{ecuaciones}), (\ref{3.4v}) and (\ref{3.10viii}), we
prove that $d_\ast{X_0}=0$, that is, $X_0$ is a 1-cocycle of
$(\frak g ^\ast ,[\, ,\, ]^{\frak g ^\ast})$.

On the other hand, since $\Lambda$ is $\sigma$-multiplicative, we
conclude that $\epsilon =-d_\ast{}_{X_0}:\frak g\to \wedge ^2\frak g$
is a 1-cocycle with respect to the representation $ad_{(\phi _0,1)}:\frak
g\times \wedge ^2\frak g\to \wedge ^2\frak g$ (see (\ref{3.10viii})
and Proposition \ref{nuestroLu}). Therefore, (\ref{condalg1}) holds.

Now, the equality $\# _\Lambda (\delta \sigma )=\tilde{X_0} -
e^{-\sigma} \bar{X}_0$ implies that $e^{-\sigma }\bar{X}_0(\sigma
)=\tilde{X}_0(\sigma )$, i.e., $e^{-\sigma}X_0(\sigma )=X_0(\sigma )$
(note that $\sigma$ is a multiplicative function). Consequently,
\begin{equation}\label{3.10ix}
\phi _0(X_0)=X_0(\sigma )=0.
\end{equation}
Using again that $\# _\Lambda (\delta \sigma )=\tilde{X_0} -
e^{-\sigma} \bar{X}_0$ and the fact that $\sigma$ is a multiplicative
function, we obtain that 
$$
0={\cal L}_{\bar{X}}(i(\delta \sigma )(\Lambda )+e^{-\sigma}\bar{X}_0)
=i(\delta\sigma )({\cal L}_{\bar{X}}\Lambda
)-e^{-\sigma}\bar{X}(\sigma )\bar{X}_0+e^{-\sigma}[\bar{X},\bar{X}_0]
$$
for $X\in \frak g$. In particular,
\begin{equation}\label{3.10x}
\begin{array}{ccl}
0&=&\Big \{ i(\delta\sigma )({\cal L}_{\bar{X}}\Lambda
)-e^{-\sigma}\bar{X}(\sigma )\bar{X}_0+e^{-\sigma}[\bar{X},\bar{X}_0]
\Big \} (\frak e)\\
&=&i(\phi _0)((\delta _\frak e\Lambda )(X))-\phi _0(X)X_0-[X_0,X]^\frak g.
\end{array}
\end{equation}
Thus, from (\ref{3.10viii}), (\ref{3.10ix}) and (\ref{3.10x}), it
follows that $i(\phi _0)(d_\ast X)+[X_0,X]^\frak g =0$.\QED
\begin{remark}
{\rm 
If $(\Lambda ,E)$ is a Jacobi structure on a connected Lie group $G$
which satisfies the hypotheses of Theorem \ref{bajada}, $(\lcf \,
,\,\rcf _{(\Lambda ,E)},\widetilde{\#}_{(\Lambda ,E)})$ is the Lie
algebroid structure on $T^\ast G\times \R$ given by (\ref{ecjacobi})
and $\alpha ,\beta \in \frak g^\ast $ then, a direct computation
shows that, 
$$\lcf (\hat{\alpha},f),(\hat{\beta},g)\rcf _{(\Lambda ,E)}(\frak e)= 
([\alpha ,\beta ]^{\frak g^\ast},0),$$
for $(\hat{\alpha},f),(\hat{\beta},g)\in \Omega ^1(G)\times C^\infty
(G,\R )$ satisfying $\hat{\alpha}(\frak e)=\alpha,\,
\hat{\beta}(\frak e)=\beta$ and $f(\frak e)=g(\frak e)=0$.
}
\end{remark}
\begin{example}\label{ejab}
{\rm 
Let $G$ be a connected simply connected abelian Lie group of
dimension $n$ and $(\Lambda ,E)$ be a Jacobi structure on $G$ such
that $\Lambda$ is a multiplicative 2-vector and $E$ is a bi-invariant
vector field. Then, $G$ is isomorphic, as a Lie group, to the dual
space $\frak g ^\ast$ of a real vector space $\frak g$ of dimension
$n$, $\Lambda$ is a linear 2-vector on $\frak g^\ast$ and there
exists $\varphi\in\frak g ^\ast$ satisfying that $E=-C_\varphi$,
$C_\varphi$ being the constant vector field on $\frak g ^\ast$
induced by $\varphi$. Thus, from (\ref{Jacobibracket}), one can
deduce that the Jacobi bracket of two linear functions on $\frak
g^\ast$ is again linear and that the Jacobi bracket of a linear
function and the constant function 1 is a constant function.
Therefore, using the results in \cite{IM0} (see Theorem 2 and Example
1 in \cite{IM0}) we conclude that $\frak g$ is a Lie algebra with Lie
bracket $[\, ,\,] ^\frak g$ and that    
\begin{equation}\label{Jaclin}
\Lambda =\Lambda _{\frak g^{\ast}}+R\wedge C_{\varphi},\quad
		    E=-C_\varphi ,
\end{equation}
where $\Lambda _{\frak g^{\ast}}$ is the Lie-Poisson structure on 
$\frak g ^\ast$, $R$ is the radial vector field on $\frak g ^\ast$ and 
$\varphi \in \frak g ^\ast$ is a 1-cocycle of $(\frak g,[\, ,\,
]^\frak g)$. The generalized Lie bialgebra associated with the
Jacobi structure $(\Lambda ,E)$ on $\frak g ^\ast$ is $((\frak g
^\ast, 0),(\frak g ,\varphi ))$ and the Lie bracket on $\frak g
^\ast$ is trivial.  

Conversely, if $(\frak g, [\, ,\, ]^\frak g)$ is a real Lie algebra
of dimension $n$, $\varphi \in \frak g^\ast$ is a 1-cocycle of
$(\frak g, [\, ,\, ]^\frak g)$ and $(\Lambda ,E)$ is the pair given
by (\ref{Jaclin}) then $\frak g^\ast$ is a connected simply connected
abelian Lie group and $(\Lambda ,E)$ is a Jacobi structure on $\frak
g^\ast$ (see Theorem 1 in \cite{IM0}). Moreover, it is clear that
$\Lambda$ is multiplicative (linear) and that $E$ is bi-invariant
(constant). 
}
\end{example}
\section{Generalized Lie bialgebras and a generalization of the
Yang-Baxter equation method}
\setcounter{equation}{0}
From (\ref{condalg1}) and (\ref{3.8''}) we deduce that if $((\frak
g,\phi _0),(\frak g^\ast ,X_0))$ is a generalized Lie bialgebra then
$d_{\ast X_0}$ is a 1-cocycle on $\frak g$ with respect to the
representation $ad_{(\phi _0,1)}:\frak g\times \wedge ^2\frak g\to
\wedge ^2\frak g$. In this Section, we will propose a method to
obtain generalized Lie bialgebras such that $d_{\ast X_0}$ is a
1-coboundary (i.e., there exists $r\in \wedge ^2\frak g$ satisfying that  
$d_{\ast X_0}X=ad_{(\phi _0,1)}(X)(r),\mbox{ for }X\in \frak g$).
It is a generalization of the well-known Yang-Baxter equation
method to obtain Lie bialgebras (see, for instance, \cite{V}). It is
clear that our method will allow us to obtain connected Lie groups
such that their corresponding Lie algebras are generalized Lie
bialgebras.  
\begin{theorem}\label{YB}
Let $(\frak g,[\, ,\, ]^\frak g)$ be a real Lie algebra of finite
dimension. Suppose that $\phi _0\in \frak g ^\ast$ is a 1-cocycle and
that $r\in \wedge ^2\frak g$ and $X_0\in \frak g$ are such that
$$[r,r]^\frak g -2X_0\wedge r \mbox{ is }ad_{(\phi
_0,1)}\mbox{-invariant,}\quad [X_0,r]^\frak g =0,$$
$$i(\phi _0)(r)-X_0 \mbox{ is }ad_{(\phi _0,0)}\mbox{-invariant.}$$
If $[\, ,\, ]^{\frak g ^\ast}$ is the bracket on $\frak g^\ast$ given by
\begin{equation}\label{corchYB}
[\alpha ,\beta ]^{\frak g ^\ast}=coad_{\# _r(\beta )}\alpha
- coad_{\# _r(\alpha )}\beta +r(\alpha ,\beta )\phi _0 +i(X_0)(\alpha 
\wedge \beta ), 
\end{equation}
for $\alpha ,\beta \in \frak g ^\ast$, where $coad
:\frak g \times \frak g ^\ast \to \frak g ^\ast$ is the coadjoint
representation of $\frak g$ over $\frak g ^\ast$, that is, $(coad\,(X)
(\alpha ))(Y)=-\alpha [X,Y]^\frak g$, for $X,Y\in \frak g$, then $(\frak g
^\ast,[\, ,\, ]^{\frak g ^\ast})$ is a Lie algebra and the pair
$((\frak g ,\phi _0),(\frak g ^\ast ,X_0))$ is a generalized Lie
bialgebra.  
\end{theorem}
\prueba Let $G$ be a connected simply connected Lie group with Lie
algebra $\frak g$. 

We define a 2-vector $\Lambda$ and a vector field $E$ on $G$ by
\begin{equation}\label{3.21'}
 \Lambda =\tilde{r}-e^{-\sigma}\bar{r},\qquad E=-\tilde{X_0},
\end{equation}
where $\sigma$ is the unique multiplicative function satisfying that
$(\delta \sigma )(\frak e)=\phi _0$.

From Example \ref{cobordes}, we have that $\Lambda$ is a
$\sigma$-multiplicative 2-vector. On the other hand,
$$
\# _\Lambda (\delta \sigma
)-\tilde{X}_0+e^{-\sigma}\bar{X}_0
=\Big ( \widetilde{i(\phi _0)(r)-X_0}\Big )-e^{-\sigma}  \Big (
\overline{i(\phi _0)(r)-X_0}\Big ) .  
$$
Therefore, since $i(\phi _0)(r)-X_0$ is $ad_{(\phi _0,0)}$-invariant,
we obtain that
\begin{equation}\label{4.23'}
\# _\Lambda (\delta \sigma )=\tilde{X}_0-e^{-\sigma}\bar{X}_0
\end{equation}
(see Example \ref{cobordes}). Note that, from this equality, we deduce
that $\tilde{X}_0(\sigma )=\bar{X}_0(\sigma )=\phi _0(X_0)=0$ (see
the proof of Theorem \ref{bajada}).

On the other hand, using (\ref{3.21'}), (\ref{4.23'}) and the
properties of the Schouten-Nijenhuis bracket, it follows that 
$$[\Lambda ,\Lambda ]-2E\wedge \Lambda =-\Bigg \{ \Big (
\widetilde{[r,r]^\frak g -2X_0\wedge r}\Big ) -e^{-2\sigma}\Big (
\overline{[r,r]^\frak g -2X_0\wedge r}\Big )\Bigg \} .$$
Thus, since $[r,r]^\frak g -2X_0\wedge r$ is $ad_{(\phi
_0,1)}$-invariant, $[\Lambda ,\Lambda ]=2E\wedge \Lambda$. Moreover, 
$$
{\cal L}_E\Lambda 
=-{\cal L}_{\tilde{X}_0}\tilde{r}-e^{-\sigma}\tilde{X_0}(\sigma
)\bar{r}+e^{-\sigma}{\cal L}_{\tilde{X}_0}\bar{r}
=0.
$$ 
Consequently, the pair $(\Lambda ,E)$ is a Jacobi structure.
Furthermore, from (\ref{3.4v}) and (\ref{3.21'}), we deduce that
the intrinsic derivative of $\Lambda$ at $\frak e$ is given by
\begin{equation}\label{4.3'}
(\delta _\frak e\Lambda )(X)=-ad_{(\phi _0,1)}(X)(r)=-[X,r]^\frak
g+\phi _0(X)r, 
\end{equation}
for $X\in \frak g$. Using this fact and Theorem \ref{bajada}, we
conclude that the bracket on $\frak g ^\ast$ given by (\ref{corchYB})
is a Lie bracket and that the pair $((\frak g ,\phi _0),(\frak g
^\ast ,X_0))$ is a generalized Lie bialgebra.\QED
\begin{remark}
{\rm
\begin{itemize}
\item[{\it i)}] Since $d_\ast{}_{X_0}X=-(\delta _\frak e\Lambda
)(X)$, for all $X\in \frak g$, we obtain that (see (\ref{4.3'})),
$d_\ast s=[s,r]^\frak g-2X_0\wedge s-i(\phi _0)(s)\wedge r$, for all
$s\in \wedge ^2\frak g$. In particular, 
\begin{equation}\label{4.3''}
d_\ast r=[r,r]^\frak g-2X_0\wedge r-i(\phi _0)r\wedge r.
\end{equation}
\item[{\it ii)}] If $X\in \frak g$, it follows that (see (\ref{corchYB}))
\begin{equation}\label{4.3'''}
[\alpha ,\beta ]^{\frak g^\ast}(X)=-[X,r]^\frak g(\alpha ,\beta
)+r(\alpha ,\beta )\phi _0(X)+\alpha (X_0)\beta (X)-\beta (X_0)\alpha
(X).
\end{equation}
\end{itemize}
}
\end{remark}
Now, using Theorem \ref{YB}, we have
\begin{corollary}\label{corYB}
Let $(\frak g,[\, ,\, ]^\frak g)$ be a real Lie algebra of finite
dimension. Suppose that $\phi _0\in \frak g
^\ast$ is a 1-cocycle and that $r\in \wedge ^2\frak g$ and $X_0\in
\frak g$ are such that $i(\phi _0)(r)=X_0$ and $(r,X_0)$ is an
algebraic Jacobi structure on $\frak g$. If $[\, ,\, ]^{\frak g
^\ast}$ is the Lie bracket on $\frak g^\ast$ 
given by (\ref{corchYB}), then $(\frak g ^\ast,[\, ,\, ]^{\frak g ^\ast})$ is
a Lie algebra and the pair $((\frak g ,\phi _0),$ $(\frak g ^\ast
,X_0))$ is a generalized Lie bialgebra. Moreover, the linear map
$-\#_r:\frak g^\ast \to \frak g$ is a Lie algebra homomorphism.   
\end{corollary}
\prueba From Theorem \ref{YB} and Definition \ref{defialg} (see Appendix A),
we deduce that the pair $((\frak g ,\phi _0),(\frak g ^\ast ,X_0))$
is a generalized Lie bialgebra. On the other hand, if $\alpha, \beta
,\gamma\in \frak g^\ast$ then the equality $[r,r]^\frak g(\alpha
,\beta ,\gamma )=2(X_0\wedge r)(\alpha ,\beta ,\gamma )$ implies that
$\gamma [\#_r(\alpha ),\#_r(\beta )]^\frak g=[\alpha ,\beta ]^{\frak
g^\ast}(\#_r (\gamma ))$ and therefore
$$\#_r([\alpha ,\beta ]^{\frak g^\ast})=-[\#_r(\alpha ),\#_r(\beta
)]^\frak g.$$ 

\vspace{-.8cm}
\QED
\begin{remark}\label{remYB}
{\rm
Let $(\frak g, [\, ,\, ]^\frak g)$ be a real Lie algebra of finite
dimension. Assume that $(\Omega ,\omega )$ is an algebraic locally
conformal symplectic (l.c.s.) structure on $\frak g$ and denote by
$(r,X_0)$ the corresponding algebraic Jacobi structure on $\frak g$
(see Appendix A). Then, using Corollary \ref{corYB} and the fact that
$X_0=-\#_r(\omega )$, we deduce that the pair $((\frak g,-\omega
),(\frak g^\ast ,X_0))$ is a generalized Lie bialgebra. Furthermore,
since $\# _r:\frak g^\ast\to \frak g$ is a linear isomorphism (see
Appendix A), it follows that $\frak g^\ast$ is isomorphic, as a Lie
algebra, to $\frak g$.
}
\end{remark}
\section{Examples of generalized Lie bialgebras}\label{Ejemplos}
\setcounter{equation}{0}
First, we will give some examples of generalized Lie bialgebras which
are obtained using Theorem \ref{YB} and Corollary \ref{corYB}.
\subsection{Generalized Lie bialgebras from contact Lie
algebras}\label{glbcontacto} 
Let $(\frak g,[\, ,\, ]^\frak g)$ be a Lie algebra endowed with an
algebraic contact 1-form $\eta$ and let $X_0$ be the Reeb vector
of $\frak g$ (see Appendix A). If ${\cal Z}(\frak g)$ is the center of
$\frak g$ and $X\in {\cal Z}(\frak g)$ then it is clear that
$i(X)(d\eta )=0$. This implies that $X\in <X_0>$. Thus, ${\cal Z}(\frak
g)\subseteq <X_0>$ (see \cite{Di}). Therefore, we have two
possibilities: ${\cal Z}(\frak g )=\{ 0 \}$ or  ${\cal Z}(\frak g
)=<X_0>$.
 
If ${\cal Z}(\frak g )=<X_0>$ then Diatta \cite{Di} proved that
$\frak g$ is the central extension of a symplectic Lie algebra $(\frak
h ,[\, ,\, ]^\frak h )$ by $\R$ via the 2-cocycle $\Omega$, $\Omega$
being the algebraic symplectic structure on $\frak h$. Conversely, if
$(\frak h ,[\, ,\, ]^\frak h )$ is a symplectic Lie algebra, with
algebraic symplectic 2-form $\Omega$, and on the direct product $\frak
g =\frak h \oplus \R$ we consider the Lie bracket $[\, ,\, ]^\frak g$
given by
\begin{equation}\label{extcent}
[(X,\lambda ),(Y,\mu )]^\frak g =([X,Y]^\frak h ,-\Omega
(X,Y)),\mbox{ for }(X,\lambda ),(Y,\mu )\in \frak g ,
\end{equation}
then $\eta=(0,1)\in \frak h ^\ast\oplus \R\cong \frak g^\ast$ is an
algebraic contact 1-form on $\frak g$. Moreover, since
$X_0=(0,1)\in\frak h\oplus \R =\frak g$, we deduce that ${\cal
Z}(\frak g )=<X_0>$ (see \cite{Di}).

Now, suppose that $r$ is the algebraic Poisson 2-vector on $\frak h$
associated with the algebraic symplectic structure $\Omega$. Then,
the pair $(r,X_0)$ is the algebraic Jacobi structure on $\frak g$
associated with the contact 1-form $\eta$ (see (\ref{7.0}),
(\ref{7.01}), (\ref{7.04}) and (\ref{7.05}) in Appendix A). Thus, using
Theorem \ref{YB} and the fact that $X_0\in {\cal Z}(\frak g)$, we can
define a Lie bracket $[\, ,\,]^{\frak g^\ast}$ on $\frak g^\ast$ in
such a way that the pair $((\frak g ,0),(\frak g ^\ast
,X_0))$ is a generalized Lie bialgebra.

On the other hand, from Corollary \ref{corYB} and since $r$ is a
solution of the classical Yang-Baxter 
equation on $\frak h$, it follows that there exists a Lie bracket 
$[\, ,\, ]^{\frak h ^\ast}$ on $\frak h ^\ast$ in such a way that 
the pair $(\frak h ,\frak h ^\ast )$ is a Lie bialgebra. In fact, the
Lie algebras $(\frak h ,[\, ,\, ]^\frak h )$ and $(\frak h ^\ast ,[\,
,\, ]^{\frak h ^\ast})$ are isomorphic and, using (\ref{corchYB}), we
get that $[(\alpha ,\lambda ),(\beta ,\mu )]^{\frak g ^\ast}=([\alpha
,\beta ]^{\frak h ^\ast},0),$ for $(\alpha ,\lambda ),(\beta ,\mu
)\in \frak h^\ast \oplus \R\cong \frak g ^\ast$. Consequently, $\frak
g^\ast$ is isomorphic, as a Lie algebra, to the direct product $\frak
h \oplus \R$. 

We illustrate the precedent construction with two simple examples.

{\bf 1.-} Let $(\frak h ,[\, ,\, ]^\frak h)$ be the abelian Lie
algebra of dimension $2n$ and $\Omega $ the usual symplectic 2-form.
Then, $\frak h \oplus \R$ endowed with the Lie bracket given by
(\ref{extcent}) is just the Lie algebra $\frak h (1,n)$ of the
generalized Heisenberg group $\H (1,n)$ of the real matrices of the form
$$\left (
\begin{array}{ccccc}
	  1& x_1      &\ldots &x_n      &z\\
	  0&  1       &\ldots & 0       &x_{n+1}\\
  \vdots   & \vdots   &       & \vdots  &\vdots \\
	  0&  0       &\ldots & 1       &x_{2n}\\   
	  0&  0       &\ldots & 0       &1 
\end{array}
\right ) .
$$
Moreover, the 1-form $\eta$ is just the usual algebraic contact
1-form on $\frak h (1,n)$. In this case, the Lie algebra $\frak h
(1,n)^\ast$ is abelian.

{\bf 2.-} Let $(\frak h ,[\, ,\, ]^\frak h)$ be the nonabelian
solvable Lie algebra of dimension 2. We can find a basis $\{ e_1,e_2\}$
of $\frak h$ such that $[e_1,e_2]^\frak h =e_1$. If $\{ e^1 ,e^2 \}$
is the dual basis of $\frak h ^\ast$ then $\Omega =e^2 \wedge
e^1$ is an algebraic symplectic 2-form on $\frak h$. The
corresponding Lie algebra $\frak g$ admits a basis $\{ e_1,e_2,e_3\}$
satisfying $[e_1,e_2]^\frak g =e_1+e_3,\, [e_1,e_3]^\frak g
=[e_2,e_3]^\frak g =0.$
Thus, it is easy to prove that $\frak g$ is isomorphic to the direct
product of the Lie algebras $\frak h$ and $\R$. Therefore, in this
case, the Lie algebras $\frak g$ and $\frak g^\ast$ are isomorphic.    
\begin{remark}
{\rm
\begin{itemize}
\item[{\it i)}] $\frak h\oplus \R$ is the Lie algebra of the group
of linear automorphisms of $\R ^2$ which preserve a line. A. Diatta
proved in \cite{Di} that, in general, the group of linear
automorphisms of $\R^n$ which preserve a hyperplane is a contact Lie
group, that is, its Lie algebra admits an algebraic contact 1-form.
\item[{\it ii)}] A complete description of symplectic Lie algebras of
dimension 4 was obtained in \cite{MR} (for a detailed study of
symplectic Lie algebras, see also \cite{DM,LiM2}). Thus, one can
determine all contact Lie algebras of dimension 5 with center of
dimension 1 and from there, using Theorem \ref{YB}, to obtain
different examples of generalized Lie bialgebras.
\end{itemize}
}
\end{remark}
Now, we will give two examples of generalized Lie bialgebras $((\frak
g,\phi _0),(\frak g^\ast ,X_0))$ which come from an algebraic contact
structure on $\frak g$ but in both cases $\phi _0\neq 0$.
In the first example, $X_0\in {\cal Z}(\frak g)$. However, $X_0\notin
{\cal Z}(\frak g)$ in the second one. 

{\bf 1.-} Let $(\frak h,[\, ,\, ]^\frak h)$ be a symplectic Lie algebra
with symplectic 2-form $\Omega$. Moreover, suppose that
$\phi _0\in \frak h^\ast$ is a 1-cocycle on $\frak h$ such that
$i(\phi _0)r$ is $ad^\frak h_{(\phi _0,0)}$-invariant, that is,
\begin{equation}\label{ad-inv}
  [X,i(\phi _0)r]^\frak h=\phi _0(X)i(\phi _0)r,
\end{equation}
for $X\in \frak h$, where $r$ is the algebraic Poisson 2-vector
associated with the symplectic 2-form $\Omega$. If we consider
on $\frak g=\frak h\oplus\R$ the Lie bracket given by (\ref{extcent}),
it is easy to prove that $\phi _0$ is a 1-cocycle on $\frak g$. We also
have that $\eta =(0,1)\in \frak h^\ast \oplus \R\cong \frak g^\ast$
is an algebraic contact 1-form on $\frak g$ and that $(r,X_0)$ is the
corresponding Jacobi structure, where $X_0=(0,1)\in 
\frak h\oplus\R =\frak g$. On the other hand, if $(X,\lambda )\in 
\frak g$ then, using (\ref{extcent}) and (\ref{ad-inv}), we deduce that
$$
[(X,\lambda ),i(\phi _0)r-X_0]^\frak g=\phi _0(X,\lambda )(i(\phi _0)r-X_0).
$$
Therefore, $i(\phi _0)r-X_0$ is $ad^\frak g_{(\phi _0,0)}$-invariant.
Thus, from Theorem \ref{YB}, $((\frak g,\phi _0),(\frak g ^\ast,X_0))$
is a generalized Lie bialgebra.

A simple example of the precedent construction is the following one. Let
$(\frak h ,[\, ,\, ]^\frak h)$ be the nonabelian solvable Lie algebra
of dimension 2. We can find a basis $\{ e_1,e_2\}$ 
of $\frak h$ such that $[e_1,e_2]^\frak h =e_1$ and $\Omega
=e^2\wedge e^1$ is an algebraic symplectic 2-form on $\frak h$, where
$\{ e^1,e^2 \}$ is the dual basis of $\frak h ^\ast$. Then, $\phi
_0=-e^2$ is a 1-cocycle on $\frak h$ and it is easy to prove that
(\ref{ad-inv}) holds. As we know, the Lie algebra $(\frak g=\frak
h\oplus \R,[\, ,\, ]^\frak g)$ is isomorphic to the direct product
$\frak h\oplus \R$. Moreover, from (\ref{corchYB}), we deduce that,
in this case, $\frak g^\ast$ is the abelian Lie algebra of dimension 3.

{\bf 2.-} Let $(\frak g ,[\, ,\, ]^\frak g)$ be the solvable Lie
algebra of dimension 3 with basis $\{ e_1,e_2,e_3\}$ such that
$$[e_1,e_2]^\frak g=0,\quad [e_1,e_3]^\frak g=e_1,\quad[e_2,e_3]^\frak
g=-e_2.$$
Take $r=e_3\wedge (e_1-e_2)$ and $X_0=e_1+e_2$. It is easy to
prove that $(r,X_0)$ is an algebraic Jacobi structure on $\frak g$
which comes from an algebraic contact structure. Moreover, if $\{
e^1,e^2,e^3\}$ is the dual basis of $\frak g^\ast$
then $\phi _0=e^3$ is a 1-cocycle of $\frak g$ and $i(\phi _0)r-X_0$
is $ad^\frak g_{(\phi _0,0)}$-invariant. Therefore, from Theorem
\ref{YB}, we deduce that $((\frak g,\phi _0),(\frak g ^\ast,X_0))$ is
a generalized Lie bialgebra. The Lie bracket on $\frak g^\ast$ is
characterized by 
$$[e^1,e^2]^{\frak g^\ast}=0 ,\quad [e^1,e^3]^{\frak g
^\ast}=-e^3,\quad [e^2,e^3]^{\frak g^\ast}=e^3.$$
\subsection{Generalized Lie bialgebras from locally conformal symplectic Lie
algebras}\label{glblcs}
Suppose that $(r_{\frak h},X_0)$ is an algebraic contact structure on a Lie
algebra $(\frak h ,[\, ,\, ]^\frak h)$. If we consider on the direct
product of Lie algebras $\frak g =\frak h\oplus \R$ the 2-vector
\begin{equation}\label{5.3}
	   r=r_\frak h +e_0 \wedge X_0,
\end{equation}
where $e_0=(0,1)\in \frak h\oplus \R=\frak g$, then $(r, X_0)$ is an
algebraic l.c.s. structure and, using Remark
\ref{remYB}, $((\frak g,\phi _0),(\frak g ^\ast ,X_0))$ is a
generalized Lie bialgebra, with $\phi _0=(0,1)\in \frak h ^\ast
\oplus \R\cong \frak g ^\ast$. In addition, the Lie algebras $\frak
g$ and $\frak g^\ast$ are isomorphic (see Remark \ref{remYB}).
\begin{remark}
{\rm
If $H$ is a connected Lie group with Lie algebra $\frak h$ then the
pair $(\bar{r},\bar{X}_0)$ defines, on the direct product $G=H\times
\R$, a left invariant l.c.s. structure of the first kind in the sense
of Vaisman \cite{V0}.
}
\end{remark}
In the case when ${\cal Z}(\frak h)=<X_0>$ we have that the pair
$((\frak h,0),(\frak h^\ast ,X_0))$ is a generalized Lie bialgebra
(see Section \ref{glbcontacto}). Moreover, from (\ref{corchYB}) and
(\ref{5.3}), we deduce that the Lie bracket $[\, ,\,]^{\frak g
^\ast}$ on $\frak g^\ast$ can be described, in terms of the Lie
bracket $[\, ,\,]^{\frak h ^\ast}$ of $\frak h^\ast$, as follows
$$[(\alpha ,\lambda ),(\beta ,\mu )]^{\frak g^\ast}=([\alpha,
\beta]^{\frak h^\ast},r_\frak h(\alpha ,\beta )),$$
for $(\alpha ,\lambda ),(\beta ,\mu )\in \frak h^\ast \oplus \R\cong
\frak g^\ast$. Thus, since $r_\frak h$ is a 2-cocycle of the Lie
algebra $(\frak h^\ast ,[\, ,\,]^{\frak h^\ast})$ (see (\ref{4.3''})),
it follows that $\frak g^\ast$ is the central extension of $\frak
h^\ast$ by $\R$ via the 2-cocycle $r_\frak h$.

On the other hand, in \cite{Di}, Diatta proved that if $(\frak h
',[\, ,\, ]^{\frak h '})$ is an exact symplectic Lie algebra then one
can define on the direct product $\frak h =\frak h '\oplus \R$ a Lie
bracket in such a way that $\frak h$ is a contact Lie algebra, with
trivial center, and $\frak h'$ is a Lie subalgebra of $\frak h$.
Using this construction we can also obtain different examples of
generalized Lie bialgebras. Next, we will show an explicit example.

Let $\frak s\frak l (2,\R)$ be the Lie algebra of the special linear
group $SL(2,\R )$. Then, there exists a basis $\{ e_1,e_2,e_3\}$ of
$\frak s\frak l (2,\R)$ such that 
  $$[e_1,e_2]^{\frak s\frak l(2,\R)}=2e_2,\quad [e_1,e_3]^{\frak
   s\frak l(2,\R)}=-2e_3,  \quad[e_2,e_3]^{\frak s\frak l(2,\R)}=e_1.$$   
It is clear that $\frak s\frak l(2,\R)$ admits exact symplectic Lie
subalgebras and, therefore, we can apply Diatta's method in order to
obtain algebraic contact structures on $\frak s\frak l(2,\R)$. In
fact, if $\lambda ^1,\lambda ^2$ and $\lambda ^3$ are real numbers
satisfying the relation $(\lambda ^1)^2+4\lambda ^2\lambda ^3\neq 0$
then the pair $(r_{\frak s\frak l(2,\R)},X_0)$ given by
$$r_{\frak s\frak l(2,\R)}=\lambda ^1e_2\wedge e_3+\lambda
^2e_1\wedge e_2-\lambda ^3e_1\wedge e_3,\quad X_0=-(\lambda
^1e_1+2\lambda ^2e_2+2\lambda ^3e_3),$$  
defines an algebraic Jacobi structure on $\frak s\frak l (2,\R)$
which comes from an algebraic contact structure. Consequently, since
$\frak g\frak l (2,\R)$ (the Lie algebra of the general linear 
group $GL(2,\R)$) is isomorphic to the direct product $\frak s\frak l
(2,\R)\oplus \R$, we conclude that the pair $((\frak g\frak l (2,\R),\phi
_0),(\frak g\frak l (2,\R )^\ast ,$ $X_0))$ is a generalized Lie
bialgebra, where $\phi _0=(0,1)\in \frak s\frak l(2,\R)^\ast \oplus
\R\cong \frak g\frak l (2,\R)^\ast$.  

Finally, we remark that there exist examples of contact Lie algebras
with trivial center which do not admit symplectic Lie subalgebras. An
interesting case is $\frak s \frak u (2)$, the Lie algebra of the special
unitary group $SU(2)$. We can consider a basis $\{ e_1,e_2,e_3\}$ of
$\frak s \frak u (2)$ such that 
$$[e_1,e_2]^{\frak s\frak u (2)}=e_3,\quad [e_1,e_3]^{\frak s\frak u
(2)}=-e_2,\quad [e_2,e_3]^{\frak s\frak u (2)}=e_1.$$
Then, if $\lambda ^1,\lambda ^2$ and $\lambda ^3$ are real numbers, 
$(\lambda ^1,\lambda ^2,\lambda ^3)\neq (0,0,0)$, we have that
the pair $(r_{\frak s\frak u(2)},X_0)$ given by 
$$r_{\frak s\frak u(2)}=\lambda ^1e_2\wedge e_3-\lambda ^2e_1\wedge
e_3+\lambda ^3e_1\wedge e_2,\quad X_0=-(\lambda ^1e_1+\lambda
^2e_2+\lambda ^3e_3),$$  
defines an algebraic Jacobi structure on $\frak s\frak u (2)$ which
comes from an algebraic contact structure. Thus, since $\frak u (2)$
(the Lie algebra of the unitary group $U(2)$) is isomorphic to the
direct product $\frak s \frak u (2)\oplus \R$, we deduce that the
pair $((\frak u(2),\phi _0),(\frak u (2)^\ast, X_0))$ is a
generalized Lie bialgebra, where $\phi _0=(0,1)\in \frak s \frak
u(2)^\ast \oplus \R \cong\frak u (2)^\ast$.

We will treat again this example in Section \ref{Compacto}. 
\subsection{Other examples of generalized Lie bialgebras}
All the examples of generalized Lie bialgebras $((\frak g,\phi
_0),(\frak g^\ast ,X_0))$ considered in Sections
\ref{glbcontacto} and \ref{glblcs} have been obtained from an
algebraic Jacobi structure $(r,X_0)$ on $\frak g$. However, the
hypotheses of Theorem \ref{YB} do not necessarily imply that the pair
$(r,X_0)$ is an algebraic Jacobi structure on $\frak g$, as it is
shown in the following simple examples.  

{\bf 1.-} Let $\frak h (1,1)$ be the Lie algebra of the Heisenberg
group $\H (1,1)$. We have a basis $\{ e_1,e_2,e_3\}$ of $\frak h
(1,1)$ such that
$$[e_1,e_2]^{\frak h(1,1)}=e_3,\quad [e_1,e_3]^{\frak
h(1,1)}=[e_2,e_3]^{\frak h(1,1)}=0. $$ 
Take an arbitrary $r=\sum_{i<j}\lambda ^{ij}e_i\wedge e_j\in \wedge ^2\frak
h(1,1)$ and $X_0=e_3$. Then, $[r,r]^{\frak h (1,1)}-2X_0\wedge r$ is
$ad$-invariant and $X_0\in {\cal Z}(\frak h(1,1))$.
Therefore, the pair $((\frak h (1,1),0),(\frak h(1,1)^\ast ,X_0))$ is a
generalized Lie bialgebra. However, $(r,X_0)$ is not, in general, an
algebraic Jacobi structure on $\frak h (1,1)$. If $\{e^1,e^2,e^3\}$
is the dual basis of $\frak h (1,1)^\ast$ we deduce that
$$[e^1,e^2]^{\frak h(1,1)^\ast}=0,\quad [e^1,e^3]^{\frak
h(1,1)^\ast}=-(1+\lambda ^{12})e^1,\quad
[e^2,e^3]^{\frak h(1,1)^\ast}=-(1+\lambda ^{12})e^2. $$ 
{\bf 2.-} Let $\frak h$ be the abelian Lie algebra of dimension 3.
Take $\{ e_1,e_2,e_3\}$ a basis of $\frak h$ and let $\{ e^1 ,e^2
,e^3\}$ be the dual basis of $\frak h ^\ast$. Denote by $\Psi$ the
endomorphism of $\frak h$ given by $\Psi = \frac{1}{2}e_1\otimes
e^1+\frac{1}{2}e_2\otimes e^2 +e_3\otimes e^3$. $\Psi$ is a 1-cocycle
with respect to the adjoint representation of $\frak h$. Thus, we can
consider the representation of $\R$ on $\frak h$ given by $\R \times
\frak h\to \frak h$, $(\lambda ,X)\mapsto \lambda \Psi (X)$, and the
corresponding semi-direct product $\frak g =\frak h \times _\Psi \R$.
We can choose a basis $\{ e_1,e_2,e_3,e_4\}$ of $\frak g$ such that 
$$[e_4,e_1]^{\frak g}=\frac{1}{2}e_1,\quad [e_4,e_2]^{\frak
g}=\frac{1}{2}e_2,\quad [e_4,e_3]^{\frak g}=e_3,$$
and the other brackets are zero. Suppose that $\{ e^1,e^2,e^3,e^4 \}$
is the dual basis of $\frak g ^\ast$. If $r\in \wedge ^2\frak g$,
$X_0\in \frak g$ and $\phi _0\in \frak g^\ast$ are defined by 
$$r=e_1\wedge e_2-2e_3\wedge e_4,\quad X_0=e_3,\quad \phi _0=e^4 ,$$ 
then $r$, $X_0$ and $\phi _0$ satisfy the hypotheses of Theorem
\ref{YB}. However, $[r,r]^\frak g-2X_0\wedge r=2e_1\wedge e_2\wedge
e_3\neq 0$ and $i(\phi _0)r-X_0=e_3\neq 0$. Moreover, a direct
computation shows that, 
$$[e^3 ,e^4 ]^{\frak g ^\ast}=e^4,\quad [e^i ,e^j ]^{\frak g ^\ast}=0,$$
for $1\leq i<j\leq 4$, $(i,j)\neq (3,4)$. 

Finally, we will exhibit an example of a generalized Lie bialgebra
$((\frak g,\phi _0),(\frak g^\ast ,X_0))$ such that $\phi _0\neq 0$
and  $d_{\ast X_0}$ is not a 1-coboundary with respect to the
representation $ad_{(\phi _0,1)}:\frak g\times \wedge ^2\frak g\to
\wedge ^2\frak g$. Note that in Section \ref{Correspondencia} (see
Example \ref{ejab}), we obtained an example which satisfies this last
condition but, in that case, $\phi _0=0$. On the other hand, all the
examples of generalized Lie bialgebras that we have given in Section
\ref{Ejemplos} are such that $d_{\ast X_0}$ is a 1-coboundary.

Let $\frak g$ be the Lie algebra of dimension 4 with basis $\{
e_1,e_2,e_3,e_4\}$ satisfying  
$$ [e_4,e_1]^\frak g =e_1,\quad [e_4,e_2]^\frak g =e_2,\quad
[e_4,e_3]^\frak g =e_3$$ 
and the other brackets being zero. If $\{ e^1,e^2,e^3,e^4\}$ is the
dual basis of $\frak g^\ast$, we consider on $\frak g^\ast$ the Lie
bracket $[\, ,\, ]^{\frak g^\ast}$ characterized by
$$[e^1,e^2]^{\frak g^\ast}=e^3,\quad [e^1,e^4]^{\frak
g^\ast}=e^4,\quad [e^i,e^j]^{\frak g^\ast}=0,$$
for $1\leq i<j\leq 4$, $(i,j)\neq (1,2),(1,4)$. Then, the pair
$((\frak g,e^4),(\frak g^\ast ,e_1))$ is a generalized Lie bialgebra.
Moreover, it is easy to prove that there does not exist $r\in \wedge
^2\frak g$ such that $d_{\ast X_0}X=ad_{(\phi _0,1)}(X)(r)$, for all
$X\in \frak g$. 
\section{Compact generalized Lie bialgebras}\label{Compacto}
\setcounter{equation}{0}
Several authors have devoted special attention to the study of
compact Lie bialgebras and an important result in this direction is
the following one \cite{LW} (see also \cite{M}): every connected
compact semisimple Lie group has a nontrivial Poisson Lie group
structure. 

In this Section, we will describe the structure of a generalized Lie
bialgebra $((\frak g ,\phi _0),$ $(\frak g ^\ast ,X_0))$, $\frak g$
being a compact Lie algebra (that is, $\frak g$ is the Lie algebra of
a compact connected Lie group).
 
If $\phi _0=0$ and $X_0=0$, the pair $(\frak g,\frak g^\ast)$ is a
Lie bialgebra. Thus, we will suppose that $\phi _0\neq 0$ or $X_0\neq
0$. Note that if $\phi _0=0$ then $X_0\in{\cal Z}(\frak g)$ (see
(\ref{condalg3})). On the other hand, if $\phi _0\neq 0$ then we can
consider and $ad$-invariant scalar product $<\, ,\, >:\frak g\times
\frak g\to \R$ and the vector $\bar{Y}_0\in \frak g$ characterized by
the relation $\phi _0(X)=<X,\bar{Y}_0>$, for $X\in \frak g$. It is
clear that $\bar{Y}_0\neq 0$ and, moreover, using that $\phi _0$ is a
1-cocycle and the fact that $<\, ,\, >$ is and $ad$-invariant scalar
product, we obtain that $\bar{Y}_0\in {\cal Z}(\frak g)$ (we remark
that $\phi _0(Y_0)=1$ with $Y_0=\frac{\bar{Y}_0}{\phi
_0(\bar{Y}_0)}\in {\cal Z}(\frak g )$). 

Therefore, if $\phi _0\neq 0$ or $X_0\neq 0$, we have that $dim\,
{\cal Z}(\frak g )\geq 1$. This implies that a compact connected Lie
group $G$ with Lie algebra $\frak g$ cannot be semisimple.

Next, we will distinguish two cases:

{\em a) The case $\phi _0\neq 0$}

Let $\frak g$ be a compact Lie algebra and $\phi _0\in \frak g ^\ast$
a 1-cocycle, $\phi _0\neq 0$. If $\frak h$ is a Lie subalgebra of
$\frak g$ and $(r,X_0)$ is an algebraic l.c.s.
structure on $\frak h$ such that $i(\phi _0)(r)=X_0$
then, from Corollary \ref{corYB}, we deduce that the pair $((\frak g ,\phi
_0),(\frak g ^\ast ,X_0))$ is a generalized Lie bialgebra, where the
Lie bracket on $\frak g^\ast$ is given by (\ref{corchYB}).

Using the above construction, we can obtain some examples of
generalized Lie bialgebras $((\frak g ,\phi _0),(\frak g ^\ast
,X_0))$, with $\phi _0\neq 0$ and $\frak g$ a compact Lie algebra.
\begin{examples}\label{3.14}
{\rm
{\it i)} {\em Compact generalized Lie bialgebras of the first kind}.
Let $\frak g$ be a compact Lie algebra and $\frak h$ an abelian Lie
subalgebra of even dimension. Furthermore, assume that $r\in \wedge
^2\frak h$ is a nondegenerate 2-vector on $\frak h$ (that is, $r$
comes from an algebraic symplectic structure on $\frak h$) and that
$\phi _0\in \frak g^\ast$ is a 1-cocycle on $\frak g$ such that $\phi
_0\neq 0$ and $\phi _0\in \frak h ^\circ$, $\frak h ^\circ$ being the
annihilator of $\frak h$. Then, $((\frak g ,\phi _0),(\frak g ^\ast
,0))$ is a generalized Lie bialgebra. The pair $((\frak g ,\phi
_0),(\frak g ^\ast ,0))$ is said to be a compact generalized Lie
bialgebra of the first kind. 

{\it ii)} {\em Compact generalized Lie bialgebras of the second kind}.
Let $(\frak g ,[\, ,\,]^\frak g)$ be a compact real Lie algebra.
Suppose that $e_1,e_2\in\frak g$ are linearly independent and that
$[e_1,e_2]^\frak g=0$. We consider the 2-vector $r$ and the vector
$X_0$ on $\frak g$ defined by 
$$r=\lambda e_1\wedge e_2,\quad X_0=\lambda ^1e_1+\lambda ^2e_2,$$
with $\lambda\in \R-\{ 0\}$ and $(\lambda ^1,\lambda ^2 )\in \R ^2-\{
(0,0)\}$. It is clear that $(r,X_0)$ is an algebraic Jacobi structure
on $\frak g$ which comes from an algebraic l.c.s. structure on the
Lie subalgebra $\frak h=<e_1,e_2>$. Therefore, if $\phi _0\in \frak
g^\ast$ is a 1-cocycle of $\frak g$ such that $i(\phi _0)(r)=X_0$
(that is, $\phi _0(e_1)=\frac{\lambda ^2}{\lambda }$ and $\phi
_0(e_2)=-\frac{\lambda ^1}{\lambda }$) then $((\frak g,\phi
_0),(\frak g^\ast ,X_0))$ is a generalized Lie bialgebra. The pair
$((\frak g ,\phi _0),(\frak g ^\ast ,X_0))$ is said to be a compact
generalized Lie bialgebra of the second kind. 

{\it iii)} {\em Compact generalized Lie bialgebras of the third kind}.
Let $(\frak g ,[\, ,\,]^\frak g)$ be a nonabelian compact real Lie
algebra. Using the root 
space decomposition theorem (see, for instance, \cite{DK}) we have
that there exist $e_1,e_2,e_3\in \frak g$ satisfying 
\begin{equation}\label{3.21-}
[e_1,e_2]^\frak g=e_3,\quad [e_1,e_3]^\frak g=-e_2,\quad
[e_2,e_3]^\frak g=e_1.
\end{equation}
Now, suppose that $\phi _0\in \frak g ^\ast $ is a 1-cocycle on
$\frak g$ and that $e_4$ is a vector of $\frak g$ such that $\phi
_0(e_4)=1$, and $[e_4,e_i]^\frak g =0$, for $i=1,2,3$. Then, we
consider the 2-vector $r$ and the vector $X_0$ on $\frak g$ defined by
$$r=\lambda ^1(e_2\wedge e_3-e_4\wedge e_1)-\lambda ^2(e_1\wedge
e_3+e_4\wedge e_2)+\lambda ^3(e_1\wedge e_2-e_4\wedge e_3),$$
$$X_0=-(\lambda ^1e_1+\lambda ^2e_2+\lambda ^3e_3),$$
with $(\lambda ^1,\lambda ^2,\lambda ^3)\in \R ^3-\{ (0,0,0)\}$. A direct
computation proves that $(r,X_0)$ is an algebraic l.c.s. structure on
the Lie subalgebra $\frak h =<e_1,e_2,e_3,e_4>$ (see Section
\ref{glblcs}). Moreover, $i(\phi _0)(r)=X_0$. Thus, $((\frak g,\phi
_0),(\frak g ^\ast ,X_0))$ is a generalized Lie bialgebra. The pair
$((\frak g,\phi _0),(\frak g ^\ast ,X_0))$ is said to be a compact
generalized Lie bialgebra of the third kind. 
}
\end{examples}
\begin{remark}\label{3.15}
{\rm
\begin{itemize} 
\item[{\it i)}] In the third example the Lie algebra $\frak h$ is
isomorphic to $\frak u (2)$, $\frak u (2)$ being the Lie algebra of
the unitary group $U(2)$.
\item[{\it ii)}] Let $\frak g$ be a nonabelian compact Lie algebra
such that ${\cal Z}(\frak g)\neq \{ 0 \}$. Assume that
$e_1,e_2,e_3\in \frak g$ satisfy (\ref{3.21-}) and that $e_4$ is
an element of ${\cal Z}(\frak g)$, $e_4\neq 0$. If $<\, ,\, >:\frak
g\times \frak g\to \R$ is an $ad$-invariant scalar product on $\frak
g$ and $\bar{\phi}_0\in \frak g ^\ast$ is the 1-form given by
$\bar{\phi}_0(X)=<X,e_4>$, for $X\in\frak g$, then $\phi
_0=\frac{\bar{\phi}_0}{<e_4,e_4>}$ is a 1-cocycle on $\frak g$ and
$\phi _0 (e_4)=1$.   
\end{itemize}
}
\end{remark}
Next, we will show that Examples \ref{3.14} {\it i)}, {\it ii)} and
{\it iii)} are the only examples of generalized Lie bialgebras 
$((\frak g ,\phi _0),(\frak g ^\ast ,X_0))$, with $\phi _0\neq 0$ and
$\frak g$ a compact Lie algebra.

First, we will prove some results. 
\begin{theorem}\label{importante}
Let $((\frak g,\phi _0)\kern-0.7pt,\kern-0.7pt(\frak g ^\ast
,\kern-1ptX_0))$ be a generalized Lie bialgebra. Suppose that $\phi
_0(Y_0)=1$, with $Y_0\in {\cal Z}(\frak g )$. Then, there exists a
Lie subalgebra $\frak h$ of $\frak g$ and a 2-vector  $r\in \wedge
^2\frak h\subseteq \wedge ^2\frak g$ such that $X_0\in \frak h$ and:
\begin{itemize}
\item[{\it i)}] The pair $(r,X_0)$ defines an algebraic Jacobi
structure on $\frak g$ which comes from an algebraic l.c.s.
structure on $\frak h$. Moreover, $i(\phi _0)(r)=X_0$.
\item[{\it ii)}] The Lie bracket $[\, ,\, ]^{\frak g ^\ast}$ on
$\frak g ^\ast$ is given by (\ref{corchYB}).
\end{itemize}
\end{theorem}
\prueba
Denote by $r$ the 2-vector on $\frak g$ given by
\begin{equation}\label{3.21''}
r=-d_{\ast X_0}Y_0.
\end{equation}
Using (\ref{condalg2}), (\ref{condalg3}), (\ref{3.21''}) and the fact
that $Y_0\in{\cal Z}(\frak g)$, we have that
\begin{equation}\label{3.21'''}
i(\phi _0)(r)=X_0.
\end{equation}
Now, we will show that the pair $(r,X_0)$ is an algebraic Jacobi
structure on $\frak g$.

From (\ref{condalg1}), (\ref{3.21''}) and since $Y_0\in {\cal
Z}(\frak g)$, it follows that 
\begin{equation}\label{3.21iv}
0=d_{\ast X_0}[X,Y_0]^\frak g =-[X,r]^\frak g +\phi _0(X)r+d_{\ast X_0}X,
\end{equation}
for all $X\in \frak g$. Therefore, using (\ref{condalg2}),
(\ref{3.21iv}) and that $X_0$ is a 1-cocycle on $(\frak g ^\ast,[\, ,\,
]^{\frak g^\ast})$, we deduce that 
\begin{equation}\label{3.21v}
[X_0,r]^\frak g =0.
\end{equation}
On the other hand, using again (\ref{3.21iv}) and the properties of
the algebraic Schouten bracket $[\, ,\, ]^\frak g$, we conclude that
$[r',r]^\frak g=d_\ast r'+2X_0\wedge r'+r\wedge i(\phi _0)(r')$, for
$r'\in\wedge ^2\frak g$. Consequently (see (\ref{3.21'''})), 
\begin{equation}\label{3.21vi}
[r,r]^\frak g-2X_0\wedge r=d_\ast r+r\wedge X_0.
\end{equation}
But, from (\ref{3.21''}) and since $d_{\ast X_0}$ is a cohomology
operator, we have that $d_\ast r+r\wedge X_0=d_{\ast X_0}r=0$, which
implies that (see (\ref{3.21vi}))      
\begin{equation}\label{3.21vii}
[r,r]^\frak g-2X_0\wedge r=0.
\end{equation}
Thus, the pair $(r,X_0)$ is an algebraic Jacobi structure on $\frak
g$ and the rank of $(r,X_0)$ is even (see (\ref{3.21'''}),
(\ref{3.21v}) and (\ref{3.21vii})). Therefore, using Proposition
\ref{posibilidades} (see Appendix A), it follows that there exists a
Lie subalgebra $\frak h$ of $\frak g$ such that $r\in \wedge ^2\frak
h$, $X_0\in \frak h$ and the pair $(r,X_0)$ comes from an algebraic
l.c.s. structure on $\frak h$. 

Finally, from (\ref{4.3'''}) and (\ref{3.21iv}), we deduce that the
Lie bracket on $\frak g^\ast$ is given by (\ref{corchYB}).\QED

Now, we will describe the algebraic l.c.s. structures on a compact
Lie algebra.  
\begin{theorem}\label{thl.c.s.}
Let $\frak h$ be a compact Lie algebra of dimension $2k$, with $k\geq
1$. Suppose that $(r,X_0)$ is an algebraic Jacobi structure on $\frak
h$ which comes from an algebraic l.c.s. structure.
\begin{itemize}
\item[{\it i)}] If $X_0=0$ then $\frak h$ is the abelian Lie algebra
and $r\in \wedge ^2\frak h$ is a nondegenerate 2-vector on $\frak h$.
\item[{\it ii)}] If $X_0\neq 0$ and $k=1$ then $\frak h$ is the
abelian Lie algebra and $r\in \wedge ^2\frak h$ is an arbitrary
2-vector on $\frak h$, $r\neq 0$.
\item[{\it iii)}] If $X_0\neq 0$ and $k\geq 2$ then $k=2$,
$\frak h$ is isomorphic to $\frak u (2)$ and 
\begin{equation}\label{formula}
\begin{array}{c}
r=\lambda ^1(e_2\wedge e_3-e_4\wedge e_1)-\lambda ^2(e_1\wedge
e_3+e_4\wedge e_2)+\lambda ^3(e_1\wedge e_2-e_4\wedge e_3),\\
X_0=-(\lambda ^1e_1+\lambda ^2e_2+\lambda ^3e_3),
\end{array}
\end{equation}
where $(\lambda ^1,\lambda ^2,\lambda ^3)\in \R ^3-\{ (0,0,0)\}$ and
$\{ e_1,e_2,e_3,e_4\}$ is a basis of $\frak h$ such that $e_4\in
{\cal Z}(\frak h )$ and
\begin{equation}\label{corchetes}
[e_1,e_2]^\frak h =e_3,\quad[e_1,e_3]^\frak h
=-e_2,\quad[e_2,e_3]^\frak h =e_1.
\end{equation}
\end{itemize} 
\end{theorem}
\prueba Denote by $(\Omega ,\omega )$ the algebraic l.c.s.
structure on $\frak h$ associated with the pair $(r,X_0)$. Then, if
$d$ is the algebraic differential of $\frak h$, we have that (see
Appendix A) 
\begin{equation}\label{lcs}
\Omega ^k=\Omega\wedge\stackrel{(k}{\ldots}\wedge \Omega \neq 0,
\quad d\Omega =\omega\wedge \Omega ,\quad d\omega =0.
\end{equation}
{\it i)} If $X_0=0$, we obtain that $\omega =0$ and $\Omega$ is an algebraic
symplectic structure on $\frak h$. Thus, since $\frak h$ is a compact Lie
algebra, {\it i)} follows using the results in \cite{Chu} (see also
\cite{LiM2}). 

{\it ii)} It is trivial.

{\it iii)} Suppose that $X_0\neq 0$ and that $k\geq 2$. Then,
$\omega\neq 0$. Moreover, we can consider an $ad$-invariant scalar
product $<\, ,\, >:\frak h\times \frak h\to \R$ and the vector
$\bar{Y}_0$ of $\frak h$ characterized by the relation  
\begin{equation}\label{3.21x}
\omega (X)=<X,\bar{Y}_0>,\mbox{ for }X\in \frak h. 
\end{equation}
Using (\ref{3.21x}) and the fact that $\omega$ is a 1-cocycle, we
deduce that $\bar{Y}_0\in {\cal Z}(\frak h )$. Consequently,
\begin{equation}\label{3.21xi}
\omega (Y_0)=1,
\end{equation}
with $Y_0=\frac{\bar{Y}_0}{\omega (\bar{Y}_0)}\in {\cal Z}(\frak h )$.

On the other hand, if $\frak h'\subseteq \frak h$ is the annihilator
of the subspace generated by $\omega$, it is clear that $\frak h'$ is
a Lie subalgebra of $\frak h$. In fact, using (\ref{3.21xi}) and since
$Y_0\in {\cal Z}(\frak h)$ and $\omega$ is a 1-cocycle, it follows
that the linear map
$$\frak h\to \frak h '\oplus \R,\quad X\mapsto (X-\omega (X)Y_0,\omega (X)),$$
is a Lie algebra isomorphism. In addition, we will show that $\frak h
'$ admits an algebraic contact structure. For this purpose, we define
the 1-form $\bar{\eta}$ on $\frak h$ given by
\begin{equation}\label{3.21xii}
\bar{\eta}=-i(Y_0)(\Omega ).
\end{equation}
Using the equality $\omega =i(X_0)(\Omega )$, we have that 
\begin{equation}\label{3.21xiii}
\bar{\eta }(X_0)=1.
\end{equation}
Moreover, from (\ref{lcs}), (\ref{3.21xi}), (\ref{3.21xii}) and since
$Y_0\in {\cal Z}(\frak h )$, we deduce that 
\begin{equation}\label{6.15}
 0={\cal L}_{Y_0}\Omega =i(Y_0)(d\Omega )+d(i(Y_0)(\Omega ))=
\Omega +\omega\wedge\bar{\eta}-d\bar{\eta}.
\end{equation}
In particular (see (\ref{3.21xi}), (\ref{3.21xii}) and (\ref{3.21xiii}))
\begin{equation}\label{6.16''}
	   i(X_0)(d\bar{\eta})=i(Y_0)(d\bar{\eta})=0.
\end{equation}
Thus, the condition $\Omega ^k=\Omega \wedge
\stackrel{(k}{\ldots}\wedge \Omega \neq 0$ implies that
$\omega\wedge\bar{\eta}\wedge (d\bar{\eta})^{k-1}\neq 0.$  
Therefore, the restriction $\eta$ of $\bar{\eta}$ to $\frak h '$ is
an algebraic contact 1-form on $\frak h '$. Furthermore, if
$(r',X_0')$ is the algebraic Jacobi structure on $\frak h'$ associated
with the contact 1-form then, from relations
(\ref{3.21xi})-(\ref{6.16''}) and the results in Appendix A, we
obtain that $r'=r+Y_0\wedge X_0$ and $X_0'=X_0$. Consequently, taking
$e_4=-Y_0$ and using Proposition \ref{thcontacto} (see Appendix B),
we prove {\it iii)}.\QED  

Now, suppose that $((\frak g ,\phi _0),(\frak g ^\ast ,X_0))$ is a
generalized Lie bialgebra, with $\phi _0\neq 0$ and $\frak g$ a
compact Lie algebra. Under these conditions we showed, at the
beginning of this Section, that there exists $Y_0\in {\cal Z}(\frak
g)$ satisfying that $\phi _0(Y_0)=1$. Thus, using Theorem
\ref{importante}, we deduce that there exists a Lie subalgebra $\frak
h$ of $\frak g$ and a 2-vector $r\in \wedge ^2\frak h\subseteq \wedge
^2\frak g$ such that $X_0\in \frak h$, $(r,X_0)$ defines an algebraic
l.c.s. structure on $\frak h$, 
\begin{equation}\label{3.21xiv}
i(\phi _0)(r)=X_0
\end{equation}
and the Lie bracket on $\frak g ^\ast$ is given by (\ref{corchYB}).

If $X_0=0$ then, from (\ref{3.21xiv}) and Theorem \ref{thl.c.s.},
we have that $\frak h$ is an abelian Lie algebra and $\phi _0\in
\frak h^\circ$. If $X_0\neq 0$ and $dim\, \frak h=2$, $\frak h$ is
the abelian Lie algebra. On the other hand, if $X_0\neq 0$ and $dim\,
\frak h \geq 4$, using again Theorem
\ref{thl.c.s.}, we obtain that $\frak h$ is isomorphic to $\frak u
(2)$, that there exists a basis $\{ e_1,e_2,e_3,e_4\}$ of $\frak
h$ such that $e_4\in {\cal Z}(\frak h )$ and  $e_1,e_2$ and $e_3$ satisfy
(\ref{corchetes}) and that the pair $(r,X_0)$ is given by (\ref{formula}).
This implies that $\phi _0 (e_4)=1$ (see (\ref{3.21xiv})).

Thus, we have proved that Examples \ref{3.14} {\it i)}, {\it ii)} and
{\it iii)} are the only examples of generalized Lie bialgebras
$((\frak g ,\phi _0),(\frak g ^\ast ,X_0))$, with $\phi _0\neq 0$ and
$\frak g$ a compact Lie algebra. In other words,
\begin{theorem}
Let $((\frak g ,\phi _0 ),(\frak g ^\ast ,X_0))$ be a generalized Lie
bialgebra, with $\phi _0\neq 0$ and $\frak g$ a compact Lie algebra.
If $X_0=0$ (respectively, $X_0\neq 0$) then it is of the first kind
(respectively, the second or third kind).
\end{theorem}
{\em b) The case $\phi _0=0$}

We will describe the structure of a generalized Lie bialgebra
$((\frak g, 0),(\frak g ^\ast ,X_0))$, $\frak g$ being a compact Lie
algebra and $X_0\neq 0$. First, we will examine a suitable example.

Let $(\frak h ,\frak h ^\ast )$ be a Lie bialgebra and $\Psi$ be an
endomorphism of $\frak h$, $\Psi :\frak h\to \frak h$. Assume that
$\Psi$ is a 1-cocycle of $\frak h$ with respect to the adjoint
representation $ad ^\frak h:\frak h\times \frak h\to \frak h$ and
that $\Psi ^\ast -Id$ is a 1-cocycle of $\frak h^\ast$ with respect
to the adjoint representation $ad^{\frak h^\ast}:\frak h^\ast\times
\frak h ^\ast \to \frak h ^\ast$. Here, $\Psi ^\ast:\frak h^\ast \to
\frak h ^\ast$ is the adjoint homomorphism of $\Psi :\frak h\to \frak
h$. Denote by $\frak g=\frak h\oplus \R$ the direct product of the
Lie algebras $\frak h$ and $\R$ and consider on $\frak g^\ast\cong
\frak h^\ast \oplus \R$ the Lie bracket $[\, ,\, ]^{\frak g^\ast}$
defined by
\begin{equation}\label{semprod}
[(\alpha ,\lambda ),(\beta ,\mu)]^{\frak g^\ast}=([\alpha ,\beta
]^{\frak h^\ast}-\lambda (\Psi ^\ast -Id)(\beta )+\mu (\Psi ^\ast
-Id)(\alpha ),0),
\end{equation}
for $(\alpha ,\lambda ),(\beta ,\mu)\in \frak h^\ast\oplus \R\cong
\frak g^\ast$. Using (\ref{semprod}), that $(\frak h,\frak h^\ast )$
is a Lie bialgebra and the fact that $\Psi$ is a 1-cocycle, we deduce
that (\ref{condalg1}) holds. Thus, the pair $((\frak g,0),(\frak
g^\ast ,(0,1)))$ is a generalized Lie bialgebra. Moreover, it is
clear that if $\frak h$ is a compact Lie algebra then $\frak g$ is
also compact. 

Next, suppose that $\frak h$ is compact and semisimple and denote by
$d_{\frak h^\ast}$ the algebraic differential of $\frak h^\ast$.
Then, from Lemma \ref{Lu}, it follows that there exist $r\in \wedge
^2\frak h$ and $Z\in \frak h$ such that 
\begin{equation}\label{diferencial}
d_{\frak h^\ast}X=-[X,r]^\frak h,\quad \Psi (X)=[X,Z]^\frak h,\quad
\Psi ^\ast (\alpha )=coad ^\frak h_Z\alpha ,
\end{equation}
for $X\in \frak h$ and $\alpha \in \frak h^\ast$, where $coad ^\frak
h:\frak h\times \frak h^\ast \to \frak h^\ast$ is the coadjoint
representation. Using (\ref{diferencial}), we obtain that
$$[[X,Z]^\frak h ,r]^\frak h(\alpha ,\beta )=(\Psi ^\ast [\alpha
,\beta ]^{\frak h^\ast})(X),$$
$$[Z,[X,r]^\frak h]^\frak h(\alpha ,\beta )=(-[\Psi ^\ast (\alpha
),\beta]^{\frak h^\ast}-[\alpha ,\Psi ^\ast (\beta )]^{\frak h
^\ast})(X),$$ 
for $\alpha ,\beta \in \frak h^\ast$. Therefore, since $\Psi ^\ast
-Id$ is an adjoint 1-cocycle of $\frak h^\ast$, we deduce that
$$[[X,Z]^\frak h,r]^\frak h+[Z,[X,r]^\frak h]^\frak h=d_{\frak
h^\ast}X=-[X,r]^\frak h .$$ 
Consequently, the equality $[[X,Z]^\frak h,r]^\frak h+[Z,[X,r]^\frak
h]^\frak h=[X,[Z,r]^\frak h ]^\frak h $
implies that
\begin{equation}\label{ad-inva}
[X,[Z,r]^\frak h]^\frak h =-[X,r]^\frak h,\mbox{ for all }X\in \frak h.
\end{equation}
The compact character of $\frak h$ allows us to choose an $ad^\frak
h$-invariant scalar product $<\, ,\, >$ on $\frak h$. We will also
denote by $<\, ,\, >$ the natural extension of $<\, ,\, >$ to $\wedge
^2 \frak h$. This extension is a scalar product on $\wedge ^2\frak h$
and, in addition, it is easy to prove that $<[X,s]^\frak h ,t
>=$ $-<s,[X,t]^\frak h>$, for $X\in \frak h$ and $s,t \in \wedge ^2\frak
h$. Thus (see (\ref{ad-inva})), 
$$<[Z,r]^\frak h,[Z,r]^\frak h>=-<r,[Z,[Z,r]^\frak h]^\frak
h>=<r,[Z,r]^\frak h>=0,$$
i.e.,
\begin{equation}\label{igualdad}
[Z,r]^\frak h =0.
\end{equation}
Then, from (\ref{diferencial}), (\ref{ad-inva}) and (\ref{igualdad}),
we conclude that the Lie bracket $[\, ,\, ]^{\frak h^\ast}$ is
trivial which implies that
\begin{equation}\label{corch-triv}
[(\alpha ,\lambda ),(\beta ,\mu)]^{\frak g^\ast}=(\mu (coad _Z^\frak
h\alpha -\alpha)-\lambda (coad _Z^\frak h\beta -\beta ),0).
\end{equation}
\begin{remark}
{\rm
If $\frak h$ is not semisimple then the Lie bracket $[\, ,\, ]^{\frak
h^\ast}$ is not, in general, trivial. In fact, suppose that ${\cal
Z}(\frak h )\neq \{ 0\}$. We know that $\frak h$ is isomorphic, as a
Lie algebra, to the direct product $\frak h '\oplus {\cal Z}(\frak
h)$, where $\frak h'$ is a compact semisimple Lie subalgebra of
$\frak h$. Therefore, if $\Psi :\frak h\cong \frak h'\oplus {\cal
Z}(\frak h)\to \frak h\cong \frak h'\oplus {\cal Z}(\frak h)$ is the
projection on the subspace ${\cal Z}(\frak h)$, it follows that
$\Psi$ is an adjoint 1-cocycle of $\frak h$. Furthermore, if on
$(\frak h')^\ast$ we consider the trivial Lie bracket and on ${\cal
Z}(\frak h)^\ast$ an arbitrary (nontrivial) Lie bracket then the
direct product $(\frak h ')^\ast\oplus {\cal Z}(\frak h)^\ast\cong
\frak h^\ast$ is a Lie algebra, the pair $(\frak h,\frak h^\ast)$ is
a Lie bialgebra and the endomorphism $\Psi ^\ast -Id$ is an adjoint
1-cocycle of $\frak h^\ast$.
}
\end{remark}
Now, we prove
\begin{theorem}
Let $((\frak g,0),(\frak g^\ast ,X_0))$ be a generalized Lie
bialgebra with $X_0\neq 0$ and $\frak g$ a compact Lie algebra. Then:
\begin{itemize}
\item[{\it i)}] There exists a Lie subalgebra $\frak h$ of $\frak g$
such that $\frak g$ is isomorphic, as a Lie algebra, to the direct
product $\frak h\oplus \R$. Moreover, under the above isomorphism,
$\frak h^\ast$ is a Lie subalgebra of $\frak g^\ast$, the pair
$(\frak h,\frak h ^\ast )$ is a Lie bialgebra, $X_0=(0,1)\in \frak
h\oplus \R\cong \frak g$ and the Lie bracket $[\, ,\, ]^{\frak
g^\ast}$ on $\frak g^\ast$ is given by
\begin{equation}\label{semprod2}
[(\alpha ,\lambda ),(\beta ,\mu)]^{\frak g^\ast}=([\alpha ,\beta
]^{\frak h^\ast}-\lambda (\Psi ^\ast -Id)(\beta )+\mu (\Psi ^\ast
-Id)(\alpha ),0),
\end{equation}
where $\Psi\in End (\frak h)$ is an adjoint 1-cocycle of $\frak h$
and $\Psi ^\ast -Id$ is an adjoint 1-cocycle of $\frak h ^\ast$.
\item[{\it ii)}] If $dim\,{\cal Z}(\frak g)=1$ then the Lie bracket
$[\, ,\,]^{\frak h^\ast}$ is trivial and there exists $Z\in \frak h$
such that $\Psi (X)=[X,Z]^\frak h$, for all $X\in \frak h$.
\end{itemize}
\end{theorem}
\prueba {\it i)} From (\ref{condalg3}) it follows that $X_0\in {\cal
Z}(\frak g)$. We consider an $ad^\frak g$-invariant scalar product
$<\, ,\,>$ on $\frak g$ and the 1-form $\theta _0\in \frak g^\ast$
defined by $\theta _0(X)=<X,X_0>$, for all $X\in \frak g$. We have
that $\theta _0$ is a 1-cocycle of $\frak g$ and we can assume,
without the loss of generality, that $\theta _0(X_0)=1$. Then, using
(\ref{condalg1}) and the fact that $X_0$ is a 1-cocycle of $\frak
g^\ast$, we deduce that the Lie subalgebra $\frak h$ is the
annihilator of the subspace generated by $\theta _0$ and that the
endomorphism $\Psi :\frak h\to \frak h$ is given by $\Psi
(X)=X-i(\theta _0)(d_{\ast}X)$, where $d_{\ast}$ is
the algebraic differential of $\frak g^\ast$. 

{\it ii)} If $dim\, {\cal Z}(\frak g)=1$ then $\frak h$ is compact
and semisimple and the result follows.\QED
\appendix

\vspace{.5cm}
{\bf \Large Appendixes}

\vspace{-.5cm}
\section{Algebraic Jacobi structures}
\def\theequation{\Alph{section}.\arabic{equation}}
\setcounter{equation}{0}
In this Appendix, we will deal with an algebraic version of the
concept of Jacobi structure.
\begin{definition}\label{defialg}
Let $(\frak g ,[\, ,\, ]^\frak g )$ be a real Lie algebra of finite
dimension. An algebraic Jacobi structure on $\frak g$ is a pair
$(r,X_0 )$, with $r\in \wedge ^2\frak g$ and $X_0 \in \frak g$
satisfying  
      $$[r,r]^\frak g =2X_0\wedge r,\qquad [X_0 ,r]^\frak g =0,$$
where $[\, ,\, ]^\frak g$ is the algebraic Schouten bracket. 
\end{definition}
Note that the algebraic Poisson structures on $\frak g$ or, in other
words, the solutions of the classical Yang-Baxter equation on $\frak
g$ are just the algebraic Jacobi structures $(r,X_0)$ such that $X_0$
is zero.

On the other hand, let $(\frak g ,[\, ,\, ]^\frak g )$ be a real Lie
algebra of finite dimension and let $G$ be a connected Lie group with
Lie algebra $\frak g$. As we know
$[\bar{s},\bar{t}\,]=\overline{[s,t]^\frak g}$, for $s,t\in \wedge
^\ast \frak g$. Thus, if $r\in \wedge ^2\frak g$ and $X_0\in \frak g$
then $(r,X_0)$ is an algebraic Jacobi structure on $\frak g$ if and
only if the pair $(\bar{r},\bar{X}_0)$ is a left invariant Jacobi
structure on $G$. 
\begin{examples}
{\rm
{\it i)} Let $(\frak g,[\, ,\,]^\frak g)$ be a real Lie algebra of
odd dimension $2k+1$. We say that $\eta \in \frak g ^\ast$ is an {\em 
algebraic contact 1-form} on $\frak g$ if $\eta \wedge (d\eta
)^k=\eta\wedge d\eta\wedge \stackrel{(k}{\ldots}\wedge d\eta \neq 0,$
where $d$ is the algebraic differential on $\frak g$ (see \cite{Di}). In
such a case, $(\frak g ,\eta )$ is termed a {\em contact Lie algebra}.
If $(\frak g ,\eta )$ is a contact Lie algebra, we define $r\in\wedge
^2 \frak g$ and $X_0 \in \frak g$ as follows
\begin{equation}\label{7.0}
r(\alpha ,\beta )=d\eta (\flat _\eta ^{-1}(\alpha ),\flat _\eta
^{-1}(\beta)), \qquad X_0 =\flat _\eta ^{-1}(\eta ),
\end{equation}
for $\alpha ,\beta \in \frak g ^\ast$, where $\flat _\eta :\frak g\to
\frak g ^\ast$ is the isomorphism of vector spaces given by
\begin{equation}\label{7.01}	    
	    \flat _\eta (X)=i(X)(d\eta )+\eta (X)\eta,
\end{equation}
for $X\in \frak g$. The vector $X_0$ is the {\em Reeb vector} of
$\frak g$ and it is characterized by the relations
\begin{equation}\label{7.02}
i(X_0)(d\eta )=0,\quad \eta(X_0)=1.
\end{equation}
If $G$ is a connected Lie group with Lie algebra $\frak g$ then it is
clear that the left invariant 1-form $\bar{\eta}$ on $G$ satisfying
$\bar{\eta}(\frak e)=\eta$ is a contact 1-form. Moreover, the pair
$(\bar{r},\bar{X}_0)$ is just the Jacobi structure on $G$ associated
with $\bar{\eta}$ (see, for instance, \cite{DLM,GL,Li2}). Therefore,
we deduce that $(r,X_0)$ is an algebraic Jacobi structure on $\frak g$.

Using (\ref{7.0}), (\ref{7.01}) and (\ref{7.02}), we have that
$\#_r(\alpha )=-\flat _\eta ^{-1}(\alpha )+\alpha (X_0)X_0,$ for
$\alpha \in \frak g^\ast$. 

{\it ii)} Let $(\frak g,[\, ,\, ]^\frak g)$ be a real Lie algebra of
even dimension $2k$. An {\em algebraic locally conformal symplectic
(l.c.s.) structure} on $\frak g$ is a pair $(\Omega ,\omega )$, where
$\Omega \in \wedge ^2 \frak g ^\ast$ is nondegenerate (that is,
$\Omega ^k =\Omega \wedge\stackrel{(k}{\ldots}\wedge \Omega \neq 0$),
$\omega \in \frak g ^\ast$ is a 1-cocycle on $\frak g$ and $d\Omega
=\omega \wedge \Omega$. The 1-form $\omega$ is the {\em Lee 1-form}
of the l.c.s. structure. 

If $(\Omega ,\omega )$ is an algebraic l.c.s. structure on $\frak g$,
one can define $r\in \wedge ^2\frak g$ and $X_0 \in \frak g$ by
\begin{equation}\label{7.04}
r(\alpha ,\beta )=\Omega (\flat _\Omega ^{-1}(\alpha ),\flat
_\Omega ^{-1}(\beta)), \qquad X_0 =\flat _\Omega ^{-1}(\omega ),
\end{equation}
for $\alpha ,\beta \in \frak g ^\ast$, $\flat _\Omega :\frak
g\to \frak g ^\ast$ being the isomorphism of vector spaces given by
\begin{equation}\label{7.05}
      \flat _\Omega (X)=i(X)\Omega,
\end{equation}      
for $X\in \frak g$. If $G$ is a connected Lie group with Lie algebra
$\frak g$ then it is clear that the left invariant 2-form
$\bar{\Omega}$ defines a locally conformal symplectic structure on
$G$. Furthermore, the pair $(\bar{r},\bar{X}_0)$ is just the Jacobi
structure on $G$ associated with $\bar{\Omega}$ (see, for instance,
\cite{GL,K}). Consequently, we obtain that $(r,X_0 )$ is an algebraic
Jacobi structure on $\frak g$. 

In this case, using (\ref{7.04}) and (\ref{7.05}), it follows that 
$\#_r(\alpha )=-\flat _\Omega ^{-1}(\alpha )$, for $\alpha \in \frak
g^\ast$. In particular, $\#_r:\frak g^\ast \to \frak g$ is a linear
isomorphism. 

It is clear that a real Lie algebra $\frak g$ is symplectic in the
sense of \cite{LiM2} if and only if $\frak g$ is l.c.s. and the Lee
1-form is zero. Moreover, if $\frak g$ is a symplectic Lie algebra
then the 2-vector $r\in \wedge ^2\frak g$ given by (\ref{7.04}) is a
solution of the classical Yang-Baxter equation on $\frak g$.
}
\end{examples}
Now, we introduce the following definition.
\begin{definition}
Let $(\frak g,[\, ,\,]^\frak g)$ be a real Lie algebra of dimension
$n$ and $(r,X_0)$ be an algebraic Jacobi structure on $\frak g$. The
rank of $(r,X_0)$ is the dimension of the subspace $\# _r(\frak g
^\ast)+$ $<X_0>\subseteq \frak g$. Equivalently, the rank of $(r,X_0)$
is $2k\leq n$ (respectively, $2k+1\leq n$) if the rank of $r$ is $2k$
and $X_0\wedge r^k=X_0\wedge r\wedge \stackrel{(k}{\ldots}\wedge r=
0$ (respectively, $X_0\wedge r^k=X_0\wedge r\wedge
\stackrel{(k}{\ldots}\wedge r$ $\neq 0$).  
\end{definition}
If $G$ is a connected Lie group with Lie algebra $\frak g$ then it is
clear that the rank of an algebraic Jacobi structure $(r,X_0)$ on
$\frak g$ is just the rank of the Jacobi structure
$(\bar{r},\bar{X}_0)$ on $G$. Thus, the rank of a contact Lie algebra
(respectively, l.c.s. Lie algebra) of dimension $2k+1$ (respectively,
$2k$) is $2k+1$ (respectively, $2k$). Conversely, using some
well-known results about transitive Jacobi manifolds (see
\cite{DLM,GL,K}), one may prove that if $(r,X_0)$ is an algebraic
Jacobi structure of rank $2k+1$ (respectively, of rank $2k$) on a Lie
algebra $\frak g$ of dimension $2k+1$ (respectively, of dimension
$2k$) then the structure $(r,X_0)$ comes from an algebraic contact
structure (respectively, an algebraic l.c.s. structure) on $\frak g$.
Moreover,
\begin{proposition}\label{posibilidades}
Let $(\frak g, [\, ,\,]^\frak g)$ be a real Lie algebra of dimension
$n$ and $(r,X_0 )$ be an algebraic Jacobi structure on $\frak g$ of
rank $m\leq n$. Then, there exists an $m$-dimensional Lie
subalgebra $\frak h$ of $\frak g$ such that $r\in\wedge ^2\frak
h$, $X_0\in \frak h$, the pair $(r,X_0)$ defines an algebraic Jacobi
structure on $\frak h$ and:
\begin{itemize}
\item[{\it i)}] If $m$ is odd, the structure $(r,X_0)$ comes from an
algebraic contact structure on $\frak h$.
\item[{\it ii)}] If $m$ is even, the structure $(r,X_0)$ comes from
an algebraic l.c.s. structure on $\frak h$.
\end{itemize}
\end{proposition}
\prueba Let $G$ be a connected Lie group with Lie algebra $\frak g$
and $(\bar{r},\bar{X}_0)$ be the corresponding left invariant Jacobi
structure on $G$. Denote by ${\cal F}$ the characteristic foliation
on $G$ associated with the Jacobi structure $(\bar{r},\bar{X}_0)$,
that is (see \cite{DLM,GL,K}), for every $g\in G$, ${\cal F}_g$ is
the subspace of $T_gG$ defined by ${\cal F}_g=(\#
_{\bar{r}})_g(T^\ast _gG)+<\bar{X}_0(g)>$. It is clear that
\begin{equation}\label{7.1}
\bar{r}(g)\in \wedge ^2{\cal F}_g,\quad {\cal F}_g=(L_g)_\ast ({\cal
F}_\frak e),\quad dim\, {\cal F}_g=dim\, {\cal F}_\frak e=m,
\end{equation}
for all $g\in G$. Thus, using (\ref{7.1}), we deduce that $\frak
h={\cal F}_\frak e$.\QED    
\section{Compact contact Lie algebras}
\setcounter{equation}{0}
In \cite{Di}, Diatta proved that if $G$ is a Lie group which admits a
left invariant contact structure and a bi-invariant semi-Riemannian
metric then $G$ is semisimple and thus, from Theorem 5 in \cite{BW}, he
deduced that G is locally isomorphic to $SL(2,\R)$ or to $SU(2)$. 
Therefore, using this result, we have that if $\frak h$
is a compact Lie algebra endowed with an algebraic contact structure
then $\frak h$ is isomorphic to $\frak s\frak u (2)$. Next, we will
give a simple proof of this last assertion. Moreover, we will
describe all the algebraic contact structures on $\frak s\frak u(2)$.
\begin{proposition}\label{thcontacto}
Let $\frak h$ be a compact Lie algebra of dimension $2k+1$, with $k\geq
1$. Suppose that $(r,X_0)$ is an algebraic Jacobi structure on $\frak
h$ which comes from an algebraic contact structure. Then, $k=1$,
$\frak h$ is isomorphic to $\frak s \frak u (2)$ and 
$$r=\lambda ^1e_2\wedge e_3-\lambda ^2e_1\wedge e_3+\lambda
^3e_1\wedge e_2,\quad X_0=-(\lambda ^1e_1+\lambda ^2e_2+\lambda ^3e_3),$$
where $(\lambda ^1,\lambda ^2,\lambda ^3)\in \R ^3-\{ (0,0,0)\}$ and $\{
e_1,e_2,e_3\}$ is a basis of $\frak h$ such that
\begin{equation}\label{3.21viii}
[e_1,e_2]^\frak h=e_3,\quad [e_1,e_3]^\frak h=-e_2,\quad
[e_2,e_3]^\frak h=e_1. 
\end{equation}
\end{proposition}
\prueba Let $\eta$ be the algebraic contact 1-form on $\frak h$
associated with the algebraic Jacobi structure $(r,X_0)$ (see
Appendix A). We can consider an $ad$-invariant
scalar product $<\, ,\, >:\frak h\times \frak h\to \R$ on $\frak h$
and the vector $\bar{X}_0\in\frak h$ characterized by the relation
\begin{equation}\label{3.21ix}
\eta (X)=<X,\bar{X}_0>,\mbox{ for } X\in \frak h.
\end{equation}
If $d$ is the algebraic differential on $\frak h$ then, using
(\ref{3.21ix}) and the fact that $<\, ,\, >$ is an $ad$-invariant
scalar product, we have that $i(\bar{X}_0)(d\eta )=0$. This implies that
\begin{equation}\label{6.9'}
Ker(d\eta)=<X_0>=<\bar{X}_0>.
\end{equation}
Next, we will see that the rank of $\frak h$ is 1. Assume that there
exists $\xi \in \frak h$ such that $[\bar{X}_0,\xi 
]^\frak h=0$. From (\ref{3.21ix}), we obtain that $(i(\xi )d\eta
)(X)=-<\bar{X}_0,[\xi ,X]^\frak h >=0$, for all $X\in \frak g$. 
Thus, using (\ref{6.9'}), we deduce that $\bar{X}_0$ and $\xi$ are
linearly dependent. 
  
Therefore, since every abelian subset of $\frak h$ must be
contained in a maximal abelian subspace of $\frak h$, we conclude that
$<\bar{X}_0>$ is a maximal abelian subspace of $\frak h$. This
implies that the rank of $\frak h$ is 1 and $\frak h$ is isomorphic
to $\frak s\frak u (2)$ (see \cite{DK}, p. 168).

Consequently, we can consider a basis $\{ e_1,e_2,e_3\}$ of $\frak h$
which satisfies (\ref{3.21viii}). Then, it is easy to prove that an
arbitrary 1-form $\eta$ on $\frak h$, $\eta\neq 0$ and $\eta =\mu
_1e^1+\mu _2e^2+\mu _3e^3$ is an algebraic contact 1-form. Here, $\{
e^1,e^2,e^3\}$ denotes the dual basis of $\{ 
e_1,e_2,e_3\}$. The algebraic Jacobi structure $(r,X_0)$ associated
with $\eta$ is given by (see Appendix A) 
$$r=\lambda ^1e_2\wedge e_3-\lambda ^2e_1\wedge e_3+\lambda
^3e_1\wedge e_2,\quad X_0=-(\lambda ^1e_1+\lambda ^2e_2+\lambda ^3e_3)$$
with $\lambda ^i=-\frac{\mu _i}{(\mu _1^2+\mu _2^2+\mu _3^2)}$, for $i\in \{
1,2,3\}$. This ends the proof of the result. \QED

\vspace{.2cm}
{\small {\bf Acknowledgments.} Research partially supported by DGICYT
grant BFM 2000-0808. The authors wish to thank J.H. Lu for sending
her thesis to us and R. Ib\'a\~nez for calling our attention to the
thesis of A. Diatta. D. Iglesias wishes to thank Spanish Ministry of
Education and Culture for an FPU grant.}


\vspace{-.3cm}
\begin{thebibliography}{9}
{\small
\bibitem{BV} K.H. Bhaskara, K. Viswanath: {\em Poisson algebras
and Poisson manifolds}, Research Notes in Mathema\-tics, 174, Pitman,
London, 1988.

\vspace{-7pt}
\bibitem{BW} W.M. Boothby, H.C. Wang: On contact manifolds, {\em Ann.
of Math.}, {\bf 68} (1958), 721-734.

\vspace{-7pt}
\bibitem{Chu} B.-Y. Chu : Symplectic homogeneus spaces, {\em Trans.
A.M.S.}, {\bf 197} (1974), 145-159.

\vspace{-7pt}
\bibitem{CDW} A. Coste, P. Dazord, A. Weinstein: Groupo\"{\i}des
symplectiques, {\em Pub. D\'ep. Math. Lyon,} {\bf 2/A} (1987), 1-62.

\vspace{-7pt}
\bibitem{DM} J.M Dardi\'e, A. Medina: Double extension symplectique
d'un groupe de Lie symplectique, {\em Adv. Math.}, {\bf 117} (1996)
208-227. 

\vspace{-7pt}
\bibitem{DLM} P. Dazord, A. Lichnerowicz, Ch.M. Marle: Structure
locale des vari\'et\'es de Jacobi, {\em J. Math. pures et appl.}, {\bf
70} (1991), 101-152.

\vspace{-7pt}
\bibitem{Di} A. Diatta: G\'eom\'etrie de Poisson et de contact des
espaces homog\'enes, {\em Ph.D. Thesis}, University of Montpellier II
(2000).

\vspace{-7pt}
\bibitem{D} V.G. Drinfeld: Hamiltonian Lie groups, Lie bialgebras and the
geometric meaning of the classical Yang-Baxter equation, {\em Sov.
Math. Dokl.}, {\bf 27} (1983),  68-71.

\vspace{-7pt}
\bibitem{DK} J.J. Duistermaat, J.A.C. Kolk: {\em Lie Groups},
Universitext, Springer-Verlag, 2000.

\vspace{-7pt}
\bibitem{F} B. Fuchssteiner: The Lie algebra structure of degenerate
Hamiltonian and bi-Hamiltonian systems, {\em Progr. Theoret. Phys.} {\bf 68}
(1982), 1082-1104.

\vspace{-7pt}
\bibitem{GL} F. Gu\'edira, A. Lichnerowicz: G\'eom\'etrie des
alg\'ebres de Lie locales de Kirillov, {\em J. Math. pures et appl.},
{\bf 63} (1984), 407-484.

\vspace{-7pt}
\bibitem{IM0} D. Iglesias and  J.C. Marrero: Some linear Jacobi
structures on vector bundles, {\em C.R. Acad. Sci. Paris}, {\bf 331}
S\'er. I (2000), 125-130. {\em arXiv: math.DG/0007138}.

\vspace{-7pt}
\bibitem{IM} D. Iglesias and  J.C. Marrero: Generalized Lie bialgebroids 
and Jacobi structures, {\em Preprint (2000), arXiv: math.DG/0008105}.

\vspace{-7pt}
\bibitem{KS} Y. Kerbrat, Z. Souici-Benhammadi: Vari\'et\'es de Jacobi et
groupo\"{\i}des de contact, {\em C.R. Acad. Sci. Paris}, {\bf 317} S\'er. I
(1993), 81-86.

\vspace{-7pt}
\bibitem{K} A. Kirillov: Local Lie algebras, {\em Russian Math.
Surveys}, {\bf 31} (1976), 55-75.

\vspace{-7pt}
\bibitem{K-S} Y. Kosmann-Schwarzbach: Exact Gerstenhaber algebras and
Lie bialgebroids, {\em Acta Appl. Math.}, {\bf 41} (1995), 153-165.

\vspace{-7pt}
\bibitem{KM} Y. Kosmann-Schwarzbach, F. Magri: Poisson Lie groups and
complete integrability. I. Drinfeld bigebras, dual extensions and
their canonical representations, {\em Ann. Inst. H. Poincar\'e Phys.
Th\'eor.}, {\bf 49} (1988), 433-460.

\vspace{-7pt}
\bibitem{LM} P. Libermann, Ch. M. Marle: {\em Symplectic Geometry and
Analytical Mechanics}, Kluwer, Dordrecht, 1987.

\vspace{-7pt}
\bibitem{Li1} A. Lichnerowicz: Les vari\'et\'es de Poisson et leurs
alg\'ebres de Lie associ\'ees, {\em J. Differential Geometry}, {\bf 12}
(1977), 253-300.

\vspace{-7pt}
\bibitem{Li2} A. Lichnerowicz: Les vari\'et\'es de Jacobi  et leurs
alg\'ebres de Lie associ\'ees, {\em J. Math. pures et appl.}, {\bf 57}
(1978), 453-488.

\vspace{-7pt} 
\bibitem{LiM2} A. Lichnerowicz, A. Medina: On Lie groups with
left invariant symplectic or k\"ahlerian structures, {\em Letters in
Mathematical Physics}, {\bf 16} (1988), 225-235.

\vspace{-7pt}
\bibitem{L} J.-H. Lu: Multiplicative and affine Poisson structures on
Lie groups, {\em Ph.D. Thesis}, University of California at Berkeley (1990).

\vspace{-7pt}
\bibitem{LW} J.-H. Lu, A. Weinstein: Poisson Lie groups, dressing
transformations and Bruhat decompositions, {\em J. Differential
Geometry} {\bf 31} (1990), 501-526.

\vspace{-7pt}
\bibitem{Mk} K. Mackenzie: {\em Lie groupoids and Lie algebroids in
differential geometry}, Cambridge University Press, 1987.

\vspace{-7pt}
\bibitem{MX} K. Mackenzie, P. Xu: Lie bialgebroids and Poisson
groupoids, {\em Duke Math. J.} {\bf 73} (1994), 415-452.

\vspace{-7pt}
\bibitem{M} S. Majid: Matched pairs of Lie groups associated to
solutions of the Yang-Baxter equations, {\em Pacific J. of Math.}
{\bf 141} (1990), 311-332.

\vspace{-7pt}
\bibitem{MR} A. Medina, P. Revoy: Groupes de Lie \'a structure
symplectique invariante, in {\em "Symplectic Geometry, Groupoids and
Integrable Systems, S\'eminaire Sud-Rhodanien de G\'eom\'etrie"} (P.
Dazord and A. Weinstein, Eds.) Math. Sci. Res. Inst. Publ.,
Springer-Verlag, 1991, 247-266. 

\vspace{-7pt}
\bibitem{NVQ} Ng\^o-van-Qu\^e: Sur l'espace de prolongement
diff\'erentiable, {\em J. Differential Geometry}, {\bf 2} (1968), 33-40.

\vspace{-7pt}
\bibitem{Pr} J. Pradines: Th\'eorie de Lie pour les groupo\"{\i}des
diff\'erentiables. Calcul diff\'erentiel dans la cat\'egorie des
groupo\"{\i}des infinit\'esimaux, {\em C.R. Acad. Sci. Paris}, {\bf 264}
S\'er. A (1967), 245-248.  

\vspace{-7pt}
\bibitem{V0} I. Vaisman: Locally conformal symplectic manifolds, {\em
Internat. J. Math. \& Math. Sci.}, {\bf 8} (1985), 521-536. 

\vspace{-7pt}
\bibitem{V} I. Vaisman: {\em Lectures on the Geometry of Poisson
Manifolds}, Progress in Math. 118, Birkh\"auser, Basel, 1994.

\vspace{-7pt}
\bibitem{V2} I. Vaisman: The BV-algebra of a Jacobi manifold, {\em
Ann. Polon. Math.}, {\bf 73} (2000), 275-290. {\em arXiv: math.DG/
9904112}.    

\vspace{-7pt}
\bibitem{We} A. Weinstein: The local structure of Poisson manifolds,
{\em J. Differential Geometry}, {\bf 18} (1983), 523-557. Errata et
Addenda {\bf 22} (1985), 255.
}
\end{thebibliography}
\end{document}